\documentclass[12pts]{amsart}

\usepackage{amscd,latexsym,amsthm,amsfonts,amssymb,amsmath,amsxtra}
\usepackage[mathscr]{eucal}
\newcommand{\BA}{{\mathbb {A}}}

\newcommand{\BC}{{\mathbb {C}}}

\newcommand{\RX}{{\mathrm {X}}}
\newcommand{\RY}{{\mathrm {Y}}}

\newcommand{\GL}{{\mathrm{GL}}}

\newcommand{\Ind}{{\mathrm{Ind}}}

\newcommand{\I}{{\mathrm{I}}}

\newcommand{\Sp}{{\mathrm{Sp}}}

\newtheorem{thm}{Theorem}[section]
\newtheorem{cor}[thm]{Corollary}

\newtheorem{prop}[thm]{Proposition}

\begin{document}
\renewcommand{\theequation}{\arabic{equation}}
\numberwithin{equation}{section}

\title{Top Fourier coefficients of residual Eisenstein series on symplectic or metaplectic groups, induced from Speh representations}

\author{David Ginzburg}

\address{School of Mathematical Sciences, Sackler Faculty of Exact Sciences, Tel-Aviv University, Israel
69978} \email{ginzburg@tauex.tau.ac.il}

\thanks{This research was supported by the ISRAEL SCIENCE FOUNDATION
	(grant No. 461/18).}

\author{David Soudry}
\address{School of Mathematical Sciences, Sackler Faculty of Exact Sciences, Tel-Aviv University, Israel
69978} \email{soudry@tauex.tau.ac.il}





\begin{abstract}
We consider the residues at the poles in the half plane $Re(s)\geq 0$ of Eisenstein series, on symplectic groups, or their double covers, induced from Speh representations. We show that for each such pole, there is a unique maximal nilpotent orbit, attached to Fourier coefficients admitted by the corresponding residual representation. We find this orbit in each case.
\end{abstract}

\maketitle

\section{Introduction}

Let $\tau$ be an irreducible, automorphic, cuspidal representation of $\GL_n(\BA)$, where $\BA$ is the adele ring of a number field $F$. Let $\Delta(\tau,m)$ ($m$, a positive integer) denote the Speh representation of $\GL_{mn}(\BA)$, attached to $\tau$. See \cite{MW89}. This is the representation spanned by the (multi-)
residues of Eisenstein series corresponding to the parabolic induction
from
$$
\tau|\det\cdot|^{s_1}\times
\tau|\det\cdot|^{s_2}\times\cdots\times
\tau|\det\cdot|^{s_m},
$$
at the point
$$
(\frac{m-1}{2},\frac{m-3}{2},...,\frac{1-m}{2}).
$$
In this paper, we assume that $\tau$ is self-dual, and consider Eisenstein series, induced from $\Delta(\tau,m)$, on the symplectic group $\Sp_{2mn}(\BA)=\Sp^{(1)}_{2mn}$, or its double covers $\Sp^{(2)}_{2mn}(\BA)$. We will write $\Sp_{2mn}$ as a matrix group in a standard form, so that the standard Borel subgroup consists of upper triangular matrices. Let $Q_{mn}$ be the Siegel parabolic subgroup of $\Sp_{2mn}$. In the linear case, let $f_{\Delta(\tau,m),s}$ be a smooth, holomorphic section of
\begin{equation}\label{0.0.1}
\rho_{\Delta(\tau,m),s}=\Ind_{Q_{mn}(\BA)}^{\Sp_{2mn}(\BA)}\Delta(\tau,m)|\det\cdot |^s.
\end{equation}
We denote the corresponding Eisenstein series by $E(f_{\Delta(\tau,m),s})$, and sometimes also by $E^{\Sp_{2mn}}(f_{\Delta(\tau,m),s})$.
In \cite{JLZ13}, Theorem 6.2, the poles of the normalized
Eisenstein series $E^*(f_{\Delta(\tau,m),s})$ in $Re(s)\geq 0$ are determined, and they are simple. For example, when $L(\tau,\wedge^2,s)$ has a pole at $s=1$, and $L(\tau,\frac{1}{2})\neq 0$, denote this set by $\Lambda_{\tau,\wedge^2,m}=\Lambda^{(1)}_{\tau,\wedge^2,m}$. Then $\Lambda_{\tau,\wedge^2,m}$ consists of the elements
\begin{equation}\label{0.0.1'}
e_{k,m}(\wedge^2)=e^{(1)}_{k,m}(\wedge^2)=k,\  k=1,2,...,\frac{m}{2},\ \ \  m\  even,
\end{equation}
$$	
e_{k,m}(\wedge^2)=e^{(1)}_{k,m}(\wedge^2)=k-\frac{1}{2},\ k=1,2,...,\frac{m+1}{2},\ \ \ m \  odd.
$$
Note that
$$
0<e_{k,m}(\wedge^2)=\frac{m}{2},\frac{m-2}{2},\frac{m-4}{2},...
$$
In the metaplectic case, we consider the similar Eisenstein series $E(f_{\Delta(\tau,m)\gamma_\psi,s})=E^{\Sp_{2mn}^{(2)}}(f_{\Delta(\tau,m)\gamma_\psi,s})$ on $\Sp_{2mn}^{(2)}(\BA)$, corresponding to
\begin{equation}\label{0.0.2}
\rho_{\Delta(\tau,m)\gamma_\psi,s}=\Ind_{Q^{(2)}_{mn}(\BA)}^{\Sp_{2mn}^{(2)}(\BA)}\Delta(\tau,m)\gamma_\psi|\det\cdot |^s.
\end{equation}
Here, $\gamma_\psi$ is the Weil factor attached to a nontrivial character $\psi$ of $F\backslash \BA$. In the sequel, we will also denote $\gamma^{(2)}_\psi=\gamma_\psi$, $\gamma^{(1)}_\psi=1$. We will show that the analogous list of poles, which are all simple, is the following set $\Lambda^{(2)}_{\tau,\wedge^2,m}$:
\begin{equation}\label{0.0.2''}
e^{(2)}_{k,m}(\wedge^2)=k-\frac{1}{2},\  k=1,2,...,\frac{m}{2}, \ \ \ m\ even,
\end{equation}
$$
e^{(2)}_{k,m}(\wedge^2)=k,\  k=1,2,...,\frac{m-1}{2}, \ \ \ m\ odd.
$$
Thus,
$$
0<e^{(2)}_{k,m}(\wedge^2)=\frac{m-1}{2},\frac{m-3}{2},\frac{m-5}{2},...
$$
We will sometimes denote in \eqref{0.0.1'}, \eqref{0.0.2''}, $e^{(\epsilon)}_{k,m}(\wedge^2)=e^{\Sp_{2mn}^{(\epsilon)}}_{k,m}(\wedge^2)$, $\epsilon=1,2$.
There are similar lists when the symmetric square $L$-function $L(\tau,\vee^2,s)$
has a pole at $s=1$.

The main result of this paper is the existence and determination of the unique maximal nilpotent orbits attached to Fourier coefficients admitted by the residual representations at each pole in the sets $\Lambda^{(\epsilon)}_{\tau,\wedge^2,m}$, $\Lambda^{(\epsilon)}_{\tau,\vee^2,m}$.

The main tool in our proof is the second identity proved in \cite{GS18}. We review this identity in the next section. In the example above, it exhibits an Eisenstein series on $\Sp_{2in}^{(2)}(\BA)$, $E(f_{\Delta(\tau,i)\gamma_\psi,s})$, as a descent from an Eisenstein series, $E(\varphi_{\Delta(\tau,i+1),s})$, on  $\Sp_{2(i+1)n}(\BA)$. Similarly, we can express an Eisenstein series on $\Sp_{2in}(\BA)$, $E(f_{\Delta(\tau,i),s})$, as a descent from an Eisenstein series, $E(\varphi_{\Delta(\tau,i+1)\gamma_\psi,s})$, on  $\Sp^{(2)}_{2(i+1)n}(\BA)$. Now, the proof follows by induction on $i$.
Thus, we prove our main theorem simultaneusly for symplectic groups and metaplectic groups. For example, denote by   $\mathcal{E}_{\Delta(\tau,m),\wedge^2,k}$ the representation of $\Sp_{2mn}(\BA)$ generated by the residues $Res_{s=e_{k,m}(\wedge^2)}E(f_{\Delta(\tau,m),s})$ of the Eisenstein series above on $\Sp_{2mn}(\BA)$. Similarly, we consider the residual Eisenstein series $\mathcal{E}_{\Delta(\tau,m)\gamma_\psi,\wedge^2,k}$ on $\Sp^{(2)}_{2mn}(\BA)$ . Recall that partitions of $2mn$, where each odd part appears with an even multiplicity,
determine nilpotent orbits of the Lie algebra of $\Sp_{2mn}$ over the algebraic closure of $F$, and these determine Fourier coefficients along unipotent subgroups. See \cite{GRS03}. For an automorphic representation $\pi$, denote by $\mathcal{O}(\pi)$ the set of maximal partitions corresponding to ( nilpotent orbits, attached to) nontrivial Fourier coefficients admitted by $\pi$. Then\\

{\bf Theorem :} Assume that $L(\tau,\wedge^2,s)$ has a pole at $s=1$, and $L(\tau,\frac{1}{2})\neq 0$. Then, for $\epsilon=1,2$,
$$
\mathcal{O}(\mathcal{E}_{\Delta(\tau,m)\gamma_\psi^{(\epsilon)},\wedge^2,k})=((2n)^{m-2e^{(\epsilon)}_{k,m}(\wedge^2)},n^{4e^{(\epsilon)}_{k,m}(\wedge^2)}).
$$
In detail:\\
For $m$ even, $1\leq k\leq \frac{m}{2}$,
$$
\mathcal{O}(\mathcal{E}_{\Delta(\tau,m),\wedge^2,k})=((2n)^{m-2k},n^{4k}).
$$
For $m$ odd, $1\leq k\leq \frac{m+1}{2}$,
$$
\mathcal{O}(\mathcal{E}_{\Delta(\tau,m),\wedge^2,k})=((2n)^{m-2k+1},n^{4k-2}).
$$
For $m$ even, $1\leq k\leq \frac{m}{2}$,
$$
\mathcal{O}(\mathcal{E}_{\Delta(\tau,m)\gamma_\psi,\wedge^2,k})=((2n)^{m-2k+1},n^{4k-2}).
$$
For $m$ odd, $1\leq k\leq \frac{m-1}{2}$,
$$
\mathcal{O}(\mathcal{E}_{\Delta(\tau,m)\gamma_\psi,\wedge^2,k})=((2n)^{m-2k},n^{4k}).
$$
\vspace{0.1cm}

We note that the result of the theorem, when $m=1$ in the second case of the theorem,
$$
\mathcal{O}(\mathcal{E}_{\tau,\wedge^2,1})=(n^2),
$$
plays a crucial role in establishing the descent of $\tau$ to an irreducible, automorphic, cuspidal, generic representation of $\Sp^{(2)}_{2n}(\BA)$. See \cite{GRS11}. Also, the third case of the theorem, with $m=n$ and $k=1$,
$$
\mathcal{O}(\mathcal{E}_{\Delta(\tau,n)\gamma_\psi,\wedge^2,1})=((2n)^{n-1},n^2),
$$
is essential for the double descent, attaching to $\tau$ the direct sum $\oplus (\sigma\otimes \hat{\sigma})$, where $\sigma$ varies over all irreducible, automorphic, cuspidal representations which lift to $\tau$ (with respect to $\gamma_\psi$). See \cite{GS20}. We believe that the theorem will have further useful applications, certainly to automorphic descent.

We prove a similar theorem when $L(\tau,\vee^2,s)$ has a pole at $s=1$. The proof should work similarly for split orthogonal groups, but we are missing there the analog of Lemma 6 in \cite{GRS03}, which we use repeatedly. In \cite{JL16}, some of the cases of the theorem above are proved for $k$ maximal.

Since the theorem above is proved simultaneously for symplectic groups and metaplectic groups, we need to require that $L(\tau,\frac{1}{2})\neq 0$, also when we deal with $\Sp^{(2)}_{2mn}(\BA)$, when $L(\tau,\wedge^2,s)$ has a pole at $s=1$. Note, that, in this case, when we compute the constant term of $E(f_{\Delta(\tau,m)\gamma_\psi,s})$, along the Siegel radical, then the contribution of the intertwining operator, corresponding to the long Weyl element, is, up to a finite set of places $S$, containing those at infinity, and outside which the section is unramified,

\begin{equation}\label{0.0.2'}
\prod_{k=1}^{\frac{m}{2}}\frac{L^S(\tau,\wedge^2,2s-2k+2)L^S(\tau,\vee^2,2s-2k+1)}{L^S(\tau,\wedge^2,2s+2k-1)L^S(\tau,\vee^2,2s+2k)},\ \ m\ even;
\end{equation}
$$
\prod_{k=1}^{\frac{m-1}{2}}\frac{L^S(\tau,\wedge^2,2s-2k+1)}{L^S(\tau,\wedge^2,2s+2k)}\prod_{k=1}^{\frac{m+1}{2}}\frac{L^S(\tau,\vee^2,2s-2k+2)}{L^S(\tau,\vee^2,2s+2k-1)},\ \ m\ odd.
$$

The maximal possible pole of the last product, in each case (when $L(\tau,\wedge^2,s)$ has a pole at $s=1$) is at $s=\frac{m-1}{2}$. But then, when $m$ is even and $k=\frac{m}{2}$, $L^S(\tau,\vee^2,2s-2k+1)=L^S(\tau,\vee^2,2s-m+1)$ might vanish (and even to a high order) at $s=\frac{m-1}{2}$. Similarly, when $m$ is odd and $k=\frac{m+1}{2}$,
$L^S(\tau,\vee^2,2s-2k+2)=L^S(\tau,\vee^2,2s-m+1)$ might vanish at $s=\frac{m-1}{2}$.

\vspace{0.5cm}

{\bf Notation}\\
\vspace{0.1cm}

For a positive integer $k$, let $w_k$ denote the $k\times k$
permutation matrix which has $1$ along the main anti-diagonal. For a field $F'$, we write the symplectic group $\Sp_{2k}(F')$ as
$$
\Sp_{2k}(F')=\{g\in \GL_{2k} (F')\ |\
{}^tg\begin{pmatrix}&w_k\\-w_k\end{pmatrix}g=\begin{pmatrix}&w_k\\-w_k\end{pmatrix}
\}.
$$
Similarly, we have the adele group $\Sp_{2k}(\BA)$, where $\BA$ is the adele ring of the number field $F$. For a place $v$ of $F$, where $F_v\neq \BC$, we write the metaplectic group $\Sp_{2k}^{(2)}(F_v)$ according to the Ranga Rao cocycle, corresponding to the standard Siegel parabolic subgroup \cite{Rao93}. See \cite{GS18}, Sec. 1.1.

Let $r\leq k$, and let $(r_1,...,r_t)$ be a partition of $r$. We denote by $Q_{r_1,...,r_t}$ the standard parabolic subgroup of $\Sp_{2k}$, whose Levi part, $M_{r_1,...,r_t}$, is isomorphic to $\GL_{r_1}\times\cdots\times\GL_{r_t}\times \Sp_{2k-2r}$. We will denote its unipotent radical by $U_{r_1,...,r_t}$. The group $\Sp_{2k}$ will usually be clear from the context. If not, then we denote $Q^{\Sp_{2k}}_{r_1,...,r_t}$, $M^{\Sp_{2k}}_{r_1,...,r_t}$, $U^{\Sp_{2k}}_{r_1,...,r_t}$.  Similarly, in $\GL_n$, for a partition $(r_1,...,r_t)$ of $n$, we denote the corresponding standard parabolic subgroup of $\GL_n$ by $P_{r_1,...,r_t}$. We denote its Levi part and unipotent radical by $L_{r_1,...,r_t}$, $V_{r_1,...,r_t}$. We will denote, $Z_n=V_{1^n}$. This is the standard maximal unipotent subgroup of $\GL_n$.

For a matrix $a$ in $\GL_r$, $r\leq k$, we will denote
\begin{equation}\label{0.0.3}
\hat{a}=diag(a,I_{2(k-r)},a^*)\in \Sp_{2k},
\end{equation}
where $a^*=w_r{}^ta^{-1}w_r$.

We will denote the elements of the Siegel radical $U_k=U_k^{\Sp_{2k}}$ inside $\Sp_{2k}$ by
	\begin{equation}\label{0.0.3'}
	u_k(x)=\begin{pmatrix}I_k&x\\0&I_k\end{pmatrix},
	\end{equation}
	where $w_kx$ is symmetric.

We fix a nontrivial character $\psi$ of $F\backslash \BA$.
We denote by $\psi_{Z_n}$ the standard Whittaker character, corresponding to $\psi$, of $Z_n(\BA)$:
\begin{equation}\label{0.0.4}
\psi_{Z_n}(z)=\psi(z_{1,2}+z_{2,3}+\cdots+z_{n-1,n}).
\end{equation}

\section{Preliminaries and statement of the main theorem}

{\bf 1. The set of poles of the Eisenstein series $E(f_{\Delta(\tau,
		m)\gamma^{(\epsilon)}_\psi,s})$}\\
\vspace{0.1cm}

We fix a self-dual, automorphic, cuspidal representation $\tau$ of $\GL_n(\BA)$, and a positive integer $m$. Thus, either $L(\tau,\wedge^2,s)$ has a pole at $s=1$, or $L(\tau,\vee^2,s)$ has a pole at $s=1$. Note that, in the first case, $n$ must be even, and the central character of $\tau$, $\omega_\tau$, must be trivial. In the second case, $n$ can be any positive integer and $\omega_\tau^2=1$. We will assume that $\tau$ is not the trivial character of $\GL_1(\BA)$.
Let $f_{\Delta(\tau,m)\gamma^{(\epsilon)}_\psi,s}$ denote a smooth, holomorphic section of the representation \eqref{0.0.1}, or \eqref{0.0.2}, ($\epsilon=1,2$), $\rho_{\Delta(\tau,m)\gamma^{(\epsilon)}_\psi,s}$. Consider the attached Eisenstein series $E(f_{\Delta(\tau,m)\gamma^{(\epsilon)}_\psi,s})$ on $\Sp^{(\epsilon)}_{2mn}(\BA)$. In \cite{JLZ13}, Theorem 6.2, the poles of the normalized Eisenstein series on $\Sp_{2mn}(\BA)$, $E^*(f_{\Delta(\tau,m),s})$, in $Re(s)\geq 0$, are determined. We now recall this list. We denote each element of the list by $e_{k,m}(\wedge^2)$, or $e_{k,m}(\vee^2)$. We include the case of metaplectic groups which does not appear in \cite{JLZ13}. In this case, we simply form the set of poles of each partial $L$-function which appears in the numerator of \eqref{0.0.2'}. We will denote these points by $e^{(2)}_{k,m}(\wedge^2)$, or $e^{(2)}_{k,m}(\vee^2)$. We will prove later that these are all the poles of $E^{\Sp^{(2)}_{2mn}}(f_{\Delta(\tau,m)\gamma_\psi,s})$. At this stage, this is just a set of points.\\
\vspace{0.2cm}

{\bf Case $\wedge^2$:} Assume that $L(\tau,\wedge^2,s)$ has a pole at $s=1$ and  $L(\tau,\frac{1}{2})\neq 0$.

\begin{equation}\label{1.1}
	e_{k,m}(\wedge^2)=
\begin{cases}
k,\ \ \ \ \ \ \ k=1,2,...,\frac{m}{2},\ \ \ \ \ \ m\ even,\\
	k-\frac{1}{2},\ \ \ k=1,2,...,\frac{m+1}{2},\ \  m\ odd\ \ \ \ (1\leq k\leq [\frac{m+1}{2}]).
\end{cases}
\end{equation}

\begin{equation}\label{1.2}
	e^{(2)}_{k,m}(\wedge^2)=
\begin{cases}
	k-\frac{1}{2},\ \ \ \ \ k=1,2,...,\frac{m}{2},\ \ \ \ m\ even,\\
	k,\ \ \ \ \ \ \ \ \ \ k=1,2,...,\frac{m-1}{2} \ \ \ \ (1\leq k\leq [\frac{m}{2}]).
\end{cases}
\end{equation}

{\bf Case $\vee^2$:} Assume that $L(\tau,\vee^2,s)$ has a pole at $s=1$.

	\begin{equation}\label{1.3}
		e_{k,m}(\vee^2)=
	\begin{cases}
k-\frac{1}{2},\ \ \ \ k=1,2,...,\frac{m}{2},\ \ \ m\  even, \\
	k,\ \ \ \ \ \ \ \ \ k=1,2,...,\frac{m-1}{2},\ \ m\ odd. \ \ (1\leq k\leq [\frac{m}{2}])
	\end{cases}
	\end{equation}

	\begin{equation}\label{1.4}
		e^{(2)}_{k,m}(\vee^2)=
	\begin{cases}
	k,\ \ \ \ \ \ \ \ k=1,2,...,\frac{m}{2},\ \ \ \ \ \ m\ even,\\
	k-\frac{1}{2},\ \ \ k=1,2,...,\frac{m+1}{2},\ \  m\ odd\ \ (1\leq k\leq [\frac{m+1}{2}]).
	\end{cases}
	\end{equation}

Let $\eta$ be either $\wedge^2$ or $\vee^2$. In \eqref{1.1} - \eqref{1.4}, we will sometimes denote
$$
e^{(\epsilon)}_{k,m}(\eta)=e^{\Sp^{(\epsilon)}_{2mn}}_{k,m}(\eta),
$$
in order to recall the group in question. Denote by $\Lambda^{(\epsilon)}_{\tau,\eta,m}$ the set of points $e^{(\epsilon)}_{k,m}(\eta)$ listed above, in each case. Denote the normalizing factor of the Eisenstein series above, $E(f_{\Delta(\tau,m)\gamma^{(\epsilon)}_\psi,s})$, by $D_\tau^{(\epsilon)}(s)$. It is easy to check that $D_\tau^{(\epsilon)}(s)$  is holomorphic and nonzero at each $e^{(\epsilon)}_{k,m}(\eta)$, and so we may replace  $E^*(f_{\Delta(\tau,m)\gamma^{(\epsilon)}_\psi,s})$ by $E(f_{\Delta(\tau,m)\gamma^{(\epsilon)}_\psi,s})$.

\begin{prop}\label {prop (0)}
Assume that $L(\tau,\vee^2,s)$ has a pole at $s=1$. The Eisenstein $E(f_{\Delta(\tau,m)\gamma_\psi,s})$, on $\Sp^{(2)}_{2mn}(\BA)$,	
has a simple pole at $s=\frac{m}{2}$, as the section varies.
\end{prop}
\begin{proof}
The proof is straightforward, by examining the constant term along the Siegel radical $U_{mn}$, and showing that it has a pole at $s=\frac{m}{2}$.
We have, for $Re(s)$ sufficiently large, the standard expression of the constant term, along $U_{mn}$,
\begin{multline}\label{0.1}
E^{U_{mn}}(f_{\Delta(\tau,m)\gamma_\psi,s})(1)=\\
\sum_{w\in Q_{mn}(F)\backslash \Sp_{2mn}(F)/Q_{mn}(F)}\sum_{\gamma\in M^w_{mn}(F)\backslash M_{mn}(F)}\int_{U_{mn}^w(F)\backslash U_{mn}(\BA)} 	
f_{\Delta(\tau,m)\gamma_\psi,s}((wu\gamma,1))du,
\end{multline}	
where $M_{mn}^w=M_{mn}\cap w^{-1}Q_{mn}w$ and $U^w_{mn}=U_{mn}\cap w^{-1}Q_{mn}w$. Recall that $\Sp_{2mn}(F)$ and $U_{mn}(\BA)$ split in $\Sp^{(2)}_{2mn}(\BA)$.
We can choose the following representatives $w=\epsilon_r$, $0\leq r\leq mn$,
$$
\epsilon_r=\begin{pmatrix}I_r\\&0&I_{mn-r}\\&-I_{mn-r}&0\\&&&I_r\end{pmatrix}.
$$
The calculation of $M_{mn}^{\epsilon_r}$ and $U_{mn}^{\epsilon_r}$ appears in \cite{GRS11}, Sec. 4.3. Thus, $M_{mn}^{\epsilon_r}=\widehat{P_{r,mn-r}}$ (see \eqref{0.0.3}). Factoring the $du$-integration in \eqref{0.1}, for each $r$, through $U_{mn}^{\epsilon_r}(F)\backslash U_{mn}^{\epsilon_r}(\BA)$, and conjugating the elements of $U_{mn}^{\epsilon_r}(\BA)$ by $\epsilon_r$, the corresponding $du$-integral, contains the following inner integral, for given $r$, $u\in U_{mn}^{\epsilon_r}(\BA)\backslash U_{mn}(\BA)$, $\gamma\in P_{r,mn-r}(F)\backslash \GL_{mn}(F)$,
\begin{multline}\label{0.2}
f_{\Delta(\tau,m)\gamma_\psi,s}^{V_{r,mn-r}}((\epsilon_r u\hat{\gamma},1))=\\
\int_{M_{r\times (mn-r)}(F)\backslash M_{r\times (mn-r)}(\BA)}f_{\Delta(\tau,m)\gamma_\psi,s}((\hat{v}_{r,mn-r}(x)\epsilon_r u\hat{\gamma},1))dx.
\end{multline}
This is is an application of the constant term along $V_{r,mn-r}$, applied to an element of $\Delta(\tau,m)$. For this to be nonzero, $r$ must be a multiple of $n$. This follows from the cuspidality of $\tau$. Thus in \eqref{0.1}, only $r=jn$, $0\leq j\leq m$ contribute. We get
\begin{multline}\label{0.3}
E^{U_{mn}}(f_{\Delta(\tau,m),s})(1)=\\
\sum_{j=0}^m\sum_\gamma\int_{S^\eta_{(m-j)n}(\BA)} 	
f_{\Delta(\tau,m)\gamma_\psi,s}^{V_{r,mn-r}}((\epsilon_{jn}\begin{pmatrix}I_{jn}\\&u_{(m-j)n}(z)\\&&I_{jn}\end{pmatrix}\hat{\gamma},1))du.
\end{multline} 	
The inner sum is over $\gamma\in P_{jn,(m-j)n}(F)\backslash \GL_{mn}(F)$, and $S^\eta_{(m-j)n}$ consists of the $(m-j)n\times (m-j)n$ matrices $z$ satisfying ${}^t(w_{(m-j)n}z)=(w_{(m-j)n}z)$.

Assume that the section $f_{\Delta(\tau,m)\gamma_\psi,s}$ is decomposable, and fix a finite set of places $S$, containing the Archimedean places, such that outside $S$, it is unramified (as well as $\tau$). Fix a place $v_0\in S$, and assume that the local section at $v_0$ is supported inside the open cell $Q^{(2)}_{mn}(F_v)(\epsilon_0U_{mn}(F_v),1)$, where $Q^{(2)}_{mn}(F_v)$ is the inverse image of $Q_{mn}(F_v)$ inside $\Sp_{2mn}(F_v)$. For such sections, only one summand in $j$ remains in \eqref{0.3}, namely the one with $j=0$, which is the intertwining operator on $\rho_{\Delta(\tau,m)\gamma_\psi,s}$, corresponding to $\epsilon_0$,
\begin{equation}\label{0.4}
E^{U_{mn}}(f_{\Delta(\tau,m)\gamma_\psi,s})(1)=M(\epsilon_0,s)(f_{\Delta(\tau,m)\gamma_\psi,s})=
\int_{U_{mn}(\BA)} 	
f_{\Delta(\tau,m)\gamma_\psi,s}((\epsilon_0u,1))du.	
\end{equation}
Since we know how to compute this intertwining operator locally on local unramified sections $f^0_{\Delta(\tau_v,m)\gamma_{\psi_v},s}$, we can directly verify that $M(\epsilon_0,s)(f_{\Delta(\tau,m),s})$ has a pole at $s=\frac{m}{2}$, and hence the Eisenstein series $E(f_{\Delta(\tau,m)\gamma_\psi,s})$ has a pole at $s=\frac{m}{2}$. In more details, denote by $\Delta(\tau_v,m)$ the component at $v$ of $\Delta(\tau,m)$, and realize it, for example, in its (local Whittaker-Speh-Shalika) model with respect to $(V_{m^n}(F_v),\psi_{V_{m^n,v}})$ (the character $\psi_{V_{m^n,v}}$ is written right after \eqref{(2)}). See \cite{CFK18}, Theorem 3. It has a unique, unramified function $W^0_{\Delta(\tau_v,m)}$, taking the value $1$ on $I_{mn}$. Think of $f^0_{\Delta(\tau_v,m)\gamma_{\psi_v},s}$ as a function $f^0_{\Delta(\tau_v,m)\gamma_{\psi_v},s}(h,g)$ on $\Sp_{2mn}^{(2)}(F_v)\times \GL^{(2)}_{mn}(F_v)$, and assume that
$$
f^0_{\Delta(\tau_v,m)\gamma_{\psi_v},s}(1,(g,1))=W^0_{\Delta(\tau_v,m)}(g)\gamma_{\psi_v}(\det(g)).
$$
Note that since $\tau$ is self-dual,  $\Delta(\tau_v,m)$ is self-dual. We have
\begin{equation}\label{0.5}
M(\epsilon_0,s)(f^0_{\Delta(\tau_v,m)\gamma_{\psi_v},s})=a(\tau_v,s)f^0_{\Delta(\tau_v,m)\gamma_{\psi_v},-s},
\end{equation}
where
$a(\tau_v,s)$ is given by \eqref{0.0.2'}. We see immediately that $a^S(\tau,s)=\prod_{v\notin S}a(\tau_v,s)$ has a pole at $s=\frac{m}{2}$. The local intertwining operators inside $S$ can be made holomorphic and nonzero by choosing appropriate local sections at the places of $S$. All in all, $M(\epsilon_0,s)(f_{\Delta(\tau,m)\gamma_\psi,s})$ has a pole at $s=\frac{m}{2}$, for a good choice of $f_{\Delta(\tau,m)\gamma_\psi,s}$. Finally, the pole at $s=\frac{m}{2}$ of $E(f_{\Delta(\tau,m)\gamma_\psi,s})$ is simple, since the local intertwining operators, coresponding to $\epsilon_0$, at the places of $S$, are holomorphic at $s=\frac{m}{2}$. The proof of this follows in the same way as that of Prop. 3.1, 3.4 (in a special case) in \cite{JLZ13}.

\end{proof}

\noindent{\bf Remark:} The same proof works for the case where $L(\tau,\wedge^2,s)$ has a pole at $s=1$ and $L(\tau,\frac{1}{2})\neq 0$, so that the corresponding Eisenstein series has a pole at $s=\frac{m}{2}$.
Our argument fails if we want to show that, if $L(\tau,\vee^2,s)$ has a pole at $s=1$, then $E(f_{\Delta(\tau,m),s})$, on $\Sp_{2mn}(\BA)$, has a pole at $s=\frac{m-1}{2}$. Indeed, in this case the factor $a^S(\tau,s)$ as in the last proof is given by
\begin{enumerate}
	\item For $m$ even,
	$$
	a^S(\tau,s)=\frac{L^S(\tau,s+\frac{1-m}{2})}{L^S(\tau,s+\frac{1+m}{2})}\prod_{k=1}^{\frac{m}{2}}\frac{L^S(\tau,\wedge^2,2s-2k+1)L^S(\tau,\vee^2,2s-2k+2)}{L^S(\tau,\wedge^2,2s+2k)L^S(\tau,\vee^2,2s+2k-1)};
	$$
	\item For $m$ odd,
	$$
	a^S(\tau,s)=\frac{L^S(\tau,s+\frac{1-m}{2})}{L^S(\tau,s+\frac{1+m}{2})}\prod_{k=1}^{\frac{m+1}{2}}\frac{L^S(\tau,\wedge^2,2s-2k+2)}{L^S(\tau,\wedge^2,2s+2k-1)}\prod_{k=1}^{\frac{m-1}{2}}\frac{L^S(\tau,\vee^2,2s-2k+1)}{L^S(\tau,\vee^2,2s+2k)};
	$$
	
\end{enumerate}
 For example, when $m$ is even, the factor $L^S(\tau,\vee^2,2s-2k+2)$, when $k=\frac{m}{2}$, has a pole at $s=\frac{m-1}{2}$. All other $L$-functions appearing in the numerator don't cancel this pole except, maybe, $L^S(\tau,s+\frac{1-m}{2})L^S(\tau,\wedge^2,2s-m+1)$, which might vanish, and even to a high order. \\
\vspace{0.1cm}

Our main theorem is on the top nilpotent orbits of the residual Eisenstein series above at the various poles listed before. The proof will also show that
the points $e^{(2)}_{k,m}(\eta)$ are indeed simple poles of the Eisenstein series $E(f_{\Delta(\tau,m)\gamma_\psi,s})$ on $\Sp_{2mn}^{(2)}(\BA)$ and there are no other poles in $Re(s)\geq 0$.
Let $\mathcal{E}_{\Delta(\tau,m)\gamma^{(\epsilon)}_\psi,\eta,k}$ be the residual representation of $\Sp^{(\epsilon)}_{2mn}(\BA)$, generated by the residues $Res_{s=e^{(\epsilon)}_{k,m}(\eta)}E(f_{\Delta(\tau,m)\gamma^{(\epsilon)}_\psi,s})$, as the section varies, at the pole $e^{(\epsilon)}_{k,m}(\eta)\in \Lambda^{(\epsilon)}_{\tau,\eta,m}$. Let $\mathcal{O}'$ be a nilpotent orbit of the Lie algebra of $\Sp_{2mn}$ over $F$. It corresponds to a partition $\underline{P}'$ of $2mn$. Assume that  $\mathcal{E}_{\Delta(\tau,m)\gamma^{(\epsilon)}_\psi,\eta,k}$ admits a nontrivial Fourier coefficient corresponding to $\mathcal{O}'$. In \cite{GS20}, Prop. 3.1, 3.2, we bounded $\mathcal{O}'$ (or $\underline{P}'$), in many cases. Denote by $\mathcal{O}(\mathcal{E}_{\Delta(\tau,m)\gamma^{(\epsilon)}_\psi,\eta,k})$ the set of maximal nilpotent orbits $\mathcal{O}$, supporting $\mathcal{E}_{\Delta(\tau,m)\gamma^{(\epsilon)}_\psi,\eta,k}$, in the sense that $\mathcal{E}_{\Delta(\tau,m)\gamma^{(\epsilon)}_\psi,\eta,k}$ admits a nontrivial Fourier coefficient corresponding to $\mathcal{O}$. Our main theorem, which will be proved in Section 4, is

\begin{thm}\label{thm 1.2}
	I. Assume that $L(\tau,\wedge^2,s)$ has a pole at $s=1$, and $L(\tau,\frac{1}{2})\neq 0$. Then the set of points \eqref{1.1}, \eqref{1.2}, $e^{(\epsilon)}_{k,m}(\wedge^2)$, is the set of poles of $E(f_{\Delta(\tau,m)\gamma^{(\epsilon)}_\psi,s})$, as the section varies, in $Re(s)\geq 0$, they are all simple, and
	$$
	\mathcal{O}(\mathcal{E}_{\Delta(\tau,m)\gamma_\psi^{(\epsilon)},\wedge^2,k})=((2n)^{m-2e^{(\epsilon)}_{k,m}(\wedge^2)},n^{4e^{(\epsilon)}_{k,m}(\wedge^2)}).
	$$
 II. Assume that $L(\tau,\vee^2,s)$ has a pole at $s=1$, and $\omega_\tau=1$. Then the set of points \eqref{1.1}, \eqref{1.2}, $e^{(\epsilon)}_{k,m}(\vee^2)$, is the set of poles of $E(f_{\Delta(\tau,m)\gamma^{(\epsilon)}_\psi,s})$, as the section varies, in $Re(s)\geq 0$, they are all simple, and\\
 \vspace{0.1cm}

 $	\mathcal{O}(\mathcal{E}_{\Delta(\tau,m)\gamma_\psi^{(\epsilon)},\vee^2,k})=$
 $$
 \begin{cases}
 ((2n)^{m-2e^{(\epsilon)}_{k,m}(\vee^2)},n^{4e^{(\epsilon)}_{k,m}(\vee^2)}),\ \ \ \ \ \ \ \ \ \ \ \ n\ even\\
 	((2n)^{m-2e^{(\epsilon)}_{k,m}(\vee^2)},(n+1)^{2e^{(\epsilon)}_{k,m}(\vee^2)},(n-1)^{2e^{(\epsilon)}_{k,m}(\vee^2)}),\ n&odd.
 	\end{cases}
 $$
 \vspace{0.1cm}

 III. Assume that $L(\tau,\vee^2,s)$ has a pole at $s=1$, and $\omega_\tau\neq 1$. Then each point $e^{(\epsilon)}_{k,m}(\vee^2)$
 is a simple pole of $E(f_{\Delta(\tau,m)\gamma^{(\epsilon)}_\psi,s}$, as the section varies, and\\
 \vspace{0.1cm}

 $	\mathcal{O}(\mathcal{E}_{\Delta(\tau,m)\gamma_\psi^{(\epsilon)},\vee^2,k})=$
 $$
 \begin{cases}
 ((2n)^{m-2e^{(\epsilon)}_{k,m}(\vee^2)},(n+2)^{2e^{(\epsilon)}_{k,m}(\vee^2)}, (n-2)^{2e^{(\epsilon)}_{k,m}(\vee^2)}),\ \ n\ even\\
 ((2n)^{m-2e^{(\epsilon)}_{k,m}(\vee^2)},(n+1)^{2e^{(\epsilon)}_{k,m}(\vee^2)},(n-1)^{2e^{(\epsilon)}_{k,m}(\vee^2)}),\ \ n\ odd.
 \end{cases}
 $$
\end{thm}

\vspace{0.5cm}

{\bf 2. Descent of Eisenstein series induced from Speh representations}\\
\vspace{0.1cm}

Our main tool of proof is the second inentity in \cite{GS18}. This identity roughly says that an appropriate descent of an Eisenstein series induced from $\Delta(\tau,i+1)$ is an Eisenstein series induced from $\Delta(\tau,i)$. This descent goes from $\Sp^{(\epsilon)}_{2n(i+1)}(\BA)$ to $\Sp^{(3-\epsilon)}_{2ni}(\BA)$. Let us recall this identity, in detail. For this identity, $\tau$ need not be self-dual. We assume that it is unitary.

Consider the unipotent radical $U_{1^{n-1}}=U^{\Sp_{2n(i+1)}}_{1^{n-1}}$. Write its elements as
\begin{equation}\label{7.1}
u=\begin{pmatrix}z&x&y\\&I_{2ni+2}&x'\\&&z^*\end{pmatrix}\in \Sp_{2n(i+1)},\  z\in Z_{n-1}.
\end{equation}
Consider the character $\psi_{n-1}=\psi_{n-1}^{\Sp_{2mn}}$ of $U_{1^{n-1}}(\BA)$ given by
\begin{equation}\label{7.2}
\psi_{n-1}(u)=\psi_{Z_{n-1}}(z)\psi(x_{n-1}\cdot e),
\end{equation}
where $\psi_{Z_{n-1}}$ is the Whittaker character \eqref{0.0.4}, and $e$ is the column vector
$$e=\begin{pmatrix}1\\0\\ \vdots\\0\end{pmatrix}\in F^{2ni+2}.
$$

The descent is via Fourier-Jacobi coefficients. Note that the character \eqref{7.2} is stabilized by the semi-direct product of $\Sp_{2ni}(\BA)$ and $t(\mathcal{H}_{2ni+1}(\BA))$, where $\Sp_{2ni}$ is realized inside $\Sp_{2n(i+1)}$ by
$$
t(h)=diag(I_n,h,I_n), h\in \Sp_{2ni},
$$
and $\mathcal{H}_{2ni+1}$ is the Heisenberge group in $2ni+1$ variables, realized inside $\Sp_{2n(i+1)}$ by
\begin{equation}\label{7.3}
t((x,e))=diag(I_{n-1},\begin{pmatrix}1&x&e\\&I_{2ni}&x'\\&&1\end{pmatrix},I_{n-1})\in \Sp_{2n(i+1)}.
\end{equation}
We have the projection $\beta$ from $U_{1^n}=U_{1^{n-1}}\rtimes t(\mathcal{H}_{2ni+1})$ onto $\mathcal{H}_{2ni+1}$. Extend the character \eqref{7.2} to $U_{1^n}(\BA)$ by making it trivial on $t(\mathcal{H}_{2ni+1}(\BA))$. We continue to denote this extension by $\psi_{n-1}$.

Let $\omega_{\psi^{-1}}$ be the Weil representation of $\mathcal{H}_{2ni+1}(\BA)\rtimes \Sp^{(2)}_{2ni}(\BA)$, associated to $\psi^{-1}$, i.e. the elements $(0,z)$ of the center of $\mathcal{H}_{2ni+1}(\BA)$ act by multiplication by $\psi^{-1}(z)$. We let $\omega_{\psi^{-1}}$ act on the space of Schwartz-Bruhat functions $\mathcal{S}(\BA^{ni})$. For $\phi\in \mathcal{S}(\BA^{ni})$, we have the corresponding theta series $\theta^\phi_{\psi^{-1}}$, viewed as a function on $\mathcal{H}_{2ni+1}(\BA)\rtimes \Sp^{(2)}_{2ni}(\BA)$.

Let $f_{\Delta(\tau,i+1)\gamma_\psi^{(\epsilon)},s}$ be a smooth, holomorphic section of $\rho_{\Delta(\tau,i+1)\gamma_\psi^{(\epsilon)},s}$. Consider the corresponding Eisenstein series $E(f_{\Delta(\tau,i+1)\gamma_\psi^{(\epsilon)},s})$, and apply to it the following Fourier-Jacobi coefficient, \\
\\
$\mathcal{D}^\phi_{\psi,ni}(E(f_{\Delta(\tau,i+1)\gamma_\psi^{(\epsilon)},s}))(h)$
\begin{equation}\label{7.5}
=\int_{U_{1^n}(F)\backslash U_{1^n}(\BA)}E(f_{\Delta(\tau,i+1)\gamma_\psi^{(\epsilon)},s},u\tilde{t}(h))\psi_{n-1}^{-1}(u)\theta^\phi_{\psi^{-1}}(\beta(u)\tilde{h})du.
\end{equation}
Here, $h\in \Sp^{(3-\epsilon)}_{2ni}(\BA)$. When $\epsilon=1$, and $h\in \Sp^{(2)}_{2ni}(\BA)$ projects to $h'\in \Sp_{2ni}(\BA)$, $\tilde{t}(h)=t(h')$, and
$\tilde{h}=h$. When $\epsilon=2$, $h\in \Sp_{2ni}(\BA)$, $\tilde{h}$ is any element of $\Sp^{(2)}_{2ni}(\BA)$, which projects to $h$, and $\tilde{t}(h)$  projects to $t(h)$, so that  the projection of $\tilde{h}$ on the second $\pm 1$ coordinate is the same as that of $\tilde{t}(h)$. Recall that in the metaplectic case, unipotent subgroups split in the double cover. Thus, we identify $U^{\Sp_{2(i+1)n}}_{1^n}(\BA)$ as a subgroup of $\Sp^{(2)}_{2(i+1)n}(\BA)$. \\
Let
\begin{equation}\label{7.6}
\alpha_0=\begin{pmatrix}0&I_{ni}&0&0\\0&0&0&I_n\\-I_n&0&0&0\\0&0&I_{ni}&0\end{pmatrix}.
\end{equation}
Denote
\begin{equation}\label{7.11}
U'_n=\{u'_{x;y}=\begin{pmatrix}I_n&x&0&y\\&I_{ni}&0&0\\&&I_{ni}&x'\\&&&I_n\end{pmatrix}\in H\}.
\end{equation}
Let, for $g \in H(\BA)$,
\begin{equation}\label{7.14}
f^\psi_{\Delta(\tau,i+1)\gamma_\psi^{(\epsilon)},s}(g)=\int_{V_{ni,1^n}(F)\backslash V_{ni,1^n}(\BA)}f_{\Delta(\tau,i+1)\gamma_\psi^{(\epsilon)},s}(\hat{v}g)\psi_{V_{ni,1^n}}(v)dv,
\end{equation}
where $\psi_{V_{ni,1^n}}$ is the character of $V_{ni,1^n}(\BA)$ given by (see \eqref{0.0.4})
\begin{equation}\label{7.11'}
\psi_{V_{ni,1^n}}(\begin{pmatrix}I_{ni}&y\\&z\end{pmatrix})=\psi_{Z_n}(z),\ z\in Z_n(\BA).
\end{equation}
Denote, for a finite set of places $S$, containing the Archimedean places, outside which $\tau$ is unramified,
$$
d_\tau^{\Sp_{2n(2j+1)},S}(s)=L^S(\tau,s+j+1)\prod_{k=1}^{j+1}L^S(\tau,\wedge^2,2s+2k-1)\prod_{k=1}^jL^S(\tau,sym^2,2s+2k);
$$
$$
d_\tau^{\Sp_{4nj},S}(s)=L^S(\tau,s+j+\frac{1}{2})\prod_{k=1}^{j}L^S(\tau,\wedge^2,2s+2k)L^S(\tau,sym^2,2s+2k-1);
$$
$$
d_\tau^{\Sp^{(2)}_{2n(2j+1)},S}(s)=\prod_{k=1}^{j}L^S(\tau,\wedge^2,2s+2k)\prod_{k=1}^{j+1}L^S(\tau,sym^2,2s+2k-1);
$$
$$
d_\tau^{\Sp^{(2)}_{4nj},S}(s)=\prod_{k=1}^{j}L^S(\tau,\wedge^2,2s+2k-1)L^S(\tau,sym^2,2s+2k);
$$

The following theorem is proved in \cite{GS18}, Theorems 7.1, 7.4, 8.1, 8.3.

\begin{thm}\label{thm 7.1}
	For $Re(s)$ sufficiently large, $h\in \Sp^{(3-\epsilon)}_{2ni}(\BA)$,
	\begin{equation}\label{7.6.1}
	\mathcal{D}^\phi_{\psi,ni}(E(f_{\Delta(\tau,i+1)\gamma_\psi^{(\epsilon)},s}))(h)=\sum_{\gamma\in Q_{ni}(F)\backslash \Sp_{2ni}(F)}\Lambda(f_{\Delta(\tau,i+1)\gamma_\psi^{(\epsilon)},s},\phi)(\gamma h),
	\end{equation}
	where
	$$
	\Lambda(f_{\Delta(\tau,i+1)\gamma_\psi^{(\epsilon)},s},\phi)((h)=\int_{U'_n(\BA)}\omega_{\psi^{-1}}(\beta(u)h)\phi(0)f^\psi_{\Delta(\tau,i+1)\gamma_\psi^{(\epsilon)},s}(\alpha_0 u\tilde{t}(h))du.
	$$
	In the sum \eqref{7.6.1}, $Q_{ni}=Q_{ni}^{\Sp_{2ni}}$. The function $\Lambda(f_{\Delta(\tau,i+1),s},\phi)$, defined for $Re(s)$ sufficiently large, by the last integral, admits analytic continuation
	to a  meromorphic function of $s$ in the whole plane. It defines a
	smooth meromorphic section of
	$$
	\rho_{\Delta(\tau,i)\gamma^{(3-\epsilon)}_{\psi^{-1}},s}=\Ind_{Q^{(3-\epsilon)}_{ni}(\BA)}^{\Sp^{(3-\epsilon)}_{2ni}(\BA)}\Delta(\tau,i)\gamma^{(3-\epsilon)}_{\psi^{-1}}
	|\det\cdot|^s .
	$$
	Thus, $\mathcal{D}^\phi_{\psi,ni}(E(f_{\Delta(\tau,i+1)\gamma_\psi^{(\epsilon)},s}))$ is an Eisenstein series on
	$\Sp^{(3-\epsilon)}_{2ni}(\BA)$, corresponding to the section $\Lambda(f_{\Delta(\tau,i+1)\gamma_\psi^{(\epsilon)},s},\phi)$ of
	$\rho_{\Delta(\tau,i)\gamma^{(3-\epsilon)}_{\psi^{-1}},s}$.\\
	Let $S$ be a finite set of places, containing the Archimedean places, outside which $\tau$ is unramified. Denote by $E^*_S(\cdot)$ an Eisenstein series normalized outside $S$. Then
	$$
	\mathcal{D}^\phi_{\psi,ni}(E^*_S(f_{\Delta(\tau,i+1)\gamma_\psi^{(\epsilon)},s}))= E^*_S(\Lambda(d_\tau^{\Sp^{(\epsilon)}_{2n(i+1)},S}(s)f_{\Delta(\tau,i+1)\gamma_\psi^{(\epsilon)},s},\phi)).
	$$
	That is, $\mathcal{D}^\phi_{\psi,ni}(E^*_S(f_{\Delta(\tau,i+1)\gamma_\psi^{(\epsilon)},s}))$, is the normalized (outside $S$) Eisenstein series on
	$\Sp^{(3-\epsilon)}_{2ni}(\BA)$ corresponding to the section $\Lambda(d_\tau^{\Sp^{(\epsilon)}_{2n(i+1)},S}(s)f_{\Delta(\tau,i+1)\gamma_\psi^{(\epsilon)},s},\phi)$.
	\end{thm}
\vspace{0.1cm}

{\bf 3. The sections $\Lambda$}\\

The sections $\Lambda(f_{\Delta(\tau,i+1)\gamma_\psi^{(\epsilon)},s},\phi)$ are decomposable, for decomposable data. Let us recall this. We will use the notation of \cite{GS18}. Fix an isomorphism
\begin{equation}\label {(1)}
p_{\tau,i}:\otimes'_v\Delta(\tau_v,i)\rightarrow \Delta(\tau,i).
\end{equation}
In the case of the double cover of $\GL_{ni}(\BA)$, we have the corresponding isomorphism (see (4.14) in \cite{GS18})
\begin{equation}\label {(2)}
\tilde{p}_{\tau,i}:\otimes'_v\Delta(\tau_v,i)\gamma_{\psi_v}\rightarrow \Delta(\tau,i)\gamma_\psi.
\end{equation}
We realize $\Delta(\tau_v,i)$ in the Whittaker-Speh-Shalika model corresponding to $\psi^{-1}_{V_{i^n},v}$ (see \cite{CFGK17}, Sec. 2.2). We will call it, for short, the $\psi^{-1}_{V_{i^n},v}$-model of $\Delta(\tau_v,i)$, and denote it by $W_{\psi^{-1}_{V_{i^n},v}}(\Delta(\tau_v,i))$. It is obtained as follows. Fix a space $E_v$ where $\Delta(\tau_v,i)$ acts. Then, up to scalar multiples, there is a unique (continuous)  linear functional $c_v$ on $E_v$, satisfying, for all $e\in E_v$, and all $u\in V_{i^n}(F_v)$,
\begin{equation}\label{(3)}
c_v(\Delta(\tau_v,i)(u)e)=\psi_{V_{i^n},v}^{-1}(u)c_v(e),
\end{equation}
where $\psi_{V_{i^n},v}$ is the character of $V_{i^n}(F_v)$ given by
$$
\psi_{V_{i^n},v}(\begin{pmatrix}I_i&x_{1,2}&x_{1,3}&\cdots&x_{1,n}\\&I_i&x_{2,3}&\cdots&x_{2,n}\\
&&\ddots\\&&&I_i&x_{n-1,n}\\&&&&I_i\end{pmatrix})=\psi_v(tr(x_{1,2}+x_{2,3}+\cdots+x_{n-1,n})).
$$
The corresponding model of $\Delta(\tau_v,i)$  is the space of functions on $\GL_{in}(F_v)$ given by $g\mapsto c_v(\Delta(\tau_v,i)(g)e)$, for $e\in E_v$. Note that when $i=1$, we get the $\psi_v^{-1}$- Whittaker model of $\tau_v$.

For each place $v$, let
\begin{equation}\label{(4)}
\rho_{\Delta(\tau_v,i;1)\gamma_{\psi_v}^{(\epsilon)},s}=\Ind_{Q^{(\epsilon)}_{ni,n}(F_v)}^{H(F_v)}(\Delta(\tau_v,i)\gamma_{\psi_v}^{(\epsilon)}|\det\cdot|^{s-\frac{1}{2}}\times
\tau_v\gamma_{\psi_v}^{(\epsilon)}|\det\cdot|^{s+\frac{i}{2}}).
\end{equation}

Consider a section $f_{\Delta(\tau_v,i;1)\gamma_{\psi_v}^{(\epsilon)},s}$ of $\rho_{\Delta(\tau_v,i;1)\gamma_{\psi_v}^{(\epsilon)},s}$. In the linear case, we view it as a
function on $\Sp_{2n(i+1)}(F_v)\times \GL_{ni}(F_v)\times \GL_n(F_v)$, such
that for a fixed element in $\Sp_{2n(i+1)}(F_v)$, the function in the two other
variables lies in the tensor product of  $W_{\psi^{-1}_{V_{i^n},v}}(\Delta(\tau_v,i))$ and the $\psi^{-1}_v$-Whittaker model of $\tau_v$. We
simplify notation and re-denote $f_{\Delta(\tau_v,i;1),s}(y)=f_{\Delta(\tau_v,i;1),s}(y;I_{ni},I_n)$. See \cite{GS18}, (4.11), for the metaplectic case.

For decomposable data,
\begin{equation}\label {(5)}
\Lambda(f_{\Delta(\tau,i+1)\gamma_\psi^{(\epsilon)},s},\phi)=p^{(\epsilon)}_{\tau,i}\circ (\otimes'_v \Lambda_v(f_{\Delta(\tau_v,i;1)\gamma_{\psi_v}^{(\epsilon)},s},\phi_v)),
\end{equation}
where, for $h\in \Sp^{(3-\epsilon)}_{2ni}(F_v)$, in the notation of Theorem \ref{thm 7.1},
\begin{equation}\label{(6)}
\Lambda_v(f_{\Delta(\tau_v,i;1)\gamma_{\psi_v}^{(\epsilon)},s},\phi_v)(h))=\int_{U'_n(F_v)}\omega_{\psi_v^{-1}}(\beta(u)h)\phi_v(0)f_{\Delta(\tau_v,i;1)\gamma_{\psi_v}^{(\epsilon)},s}(\alpha_0 u\tilde{t}(h))du.
\end{equation}
We recall that when $\epsilon=2$, we view $\Lambda_v(f_{\Delta(\tau_v,i;1)\gamma_{\psi_v},s},\phi_v)$ as an element of $\rho_{\Delta(\tau_v,i)\gamma_{\psi_v},s}$, via
$$
(a,\mu)\mapsto \Lambda_v(f_{\Delta(\tau,i;1)_v\gamma_{\psi_v},s},\phi_v)((\hat{a},\mu)h)=
$$
$$
\mu\gamma_{\psi_v}(\det(a))|\det(a)|_v^{s+\frac{ni+1}{2}}\int_{U'_n(F_v)}\omega_{\psi_v^{-1}}(\beta(u)h)\phi_v(0)f_{\Delta(\tau,i;1)_v\gamma_{\psi_v},s}(\alpha_0 u\tilde{t}(h),a;I_n)du.
$$

\begin{prop}\label{prop (1)}
For each place $v$, the local sections $\Lambda_v(f_{\Delta(\tau_v,i;1)\gamma_{\psi_v}^{(\epsilon)},s},\phi_v)$ are holomorphic in $\BC$.
\end{prop}
\begin{proof}
Consider the linear case ($\epsilon=1$). The proof follows from Prop. 7.2 in \cite{GS18}. Indeed, Eq. (7.18) in \cite{GS18} tells us that it is enough to consider the following integral
\begin{equation}\label{(14)}
\int \rho(\alpha_0)(f_{\Delta(\tau_v,i;1),s})(\begin{pmatrix}I_{ni}\\&I_n\\&y&I_n\\&&&I_{ni}\end{pmatrix})\psi_v^{-1}(y_{n,1})dy.
\end{equation}
The integration is over the lower Siegel radical of $\Sp_{2n}(F_v)$. As in the beginning of the proof of Theorem 7.3 in \cite{GS18}, the function, defined for $h\in \Sp_{2n}(F_v)$, by
$$
f_{\tau_v,s+\frac{i}{2}}(h)=\rho(\alpha_0)(f_{\Delta(\tau_v,i;1),s})(diag(I_{ni},h,I_{ni})),
$$
is a smooth, holomorphic section of $\rho_{\tau_v,s+\frac{i}{2}}$. The integral \eqref{(14)} becomes
\begin{equation}\label{(15)}
\int f_{\tau_v,s+\frac{i}{2}}(\begin{pmatrix}I_n\\y&I_n\end{pmatrix})\psi_v^{-1}(y_{n,1})dy.
\end{equation}
Clearly the integral \eqref{(15)}, which converges absolutely, for $Re(s)$ sufficiently large, is a Jacquet integral defining a Whittaker functional for
$\rho_{\tau_v,s+\frac{i}{2}}$, and hence admits an analytic continuation to a holomorphic function in $\BC$. The proof in the metaplectic case is entirely similar with simple modifications. 	
\end{proof}

\begin{cor}\label{cor (1')}
The global sections $\Lambda(f_{\Delta(\tau,i+1)\gamma_\psi^{(\epsilon)},s},\phi)$ are holomorphic at the half plane $Re(s)\geq -\frac{i}{2}$.
\end{cor}

\begin{proof}
Consider, for example, the linear case. Let $f_{\Delta(\tau,i+1),s}$ and $\phi$ be decomposable. Let $S$ be a finite set of places, containing the infinite places, outside which $f_{\Delta(\tau,i+1),s}$ and $\phi$ are unramified. At these places, assume that $\phi_v=\phi_v^0$ and $f^0_{\Delta(\tau_v,i;1),s}$ are as in Theorem 7.3 in \cite{GS18}. This theorem implies that
\begin{multline}\label{(15.1)}
\Lambda(f_{\Delta(\tau,i+1),s},\phi)=\\
\frac{1}{L^S(\tau, \wedge^2,2s+i+1)}p_{\tau,i}\circ (\otimes_{v\in S}\Lambda(f_{\Delta(\tau_v,i;1),s},\phi_v)\otimes(\otimes_{v\notin S}f^0_{\Delta(\tau_v,i)\gamma_{\psi_v}^{-1},s})).
\end{multline}
The corollary follows from Prop. \ref{prop (1)} and the fact that $L^S(\tau, \wedge^2,2s+i+1)$ doesn't vanish at $Re(2s+i+1)\geq 1$. The metaplectic case is proved similarly, using Theorem 8.2 in \cite{GS18}.
\end{proof}
\begin{prop}\label{prop (2)}
Fix a place $v$ and fix a complex number $s_0$ with $Re(s_0)\geq 0$. Consider the local sections $\Lambda_v(f_{\Delta(\tau_v,i;1)\gamma_{\psi_v}^{(\epsilon)},s_0},\phi_v)$, as in Prop. \ref{prop (1)}, at the point $s_0$. View these as maps to $\rho_{\Delta(\tau_v,i)\gamma_{\psi_v}^{(\epsilon)},s_0}$. Then, when $v$ is finite, these are surjective maps. When $v$ is Archimedean, their image is dense (in the Frechet topology).
\end{prop}
\begin{proof}
We prove the proposition in the linear case. The proof in the metaplectic case is very similar.

Denote by $\Lambda_{v,s_0}$ be the following bilinear map on $V_{\rho_{\Delta(\tau_v,i;1),s_0}}\times \mathcal{S}(F_v^{ni})$. Let $f_0$ be a function in $V_{\rho_{\Delta(\tau_v,i;1),s_0}}$, and let $\phi_v\in \mathcal{S}(F_v^{ni})$. Let $f_{\Delta(\tau_v,i;1),s}$ be any smooth, holomorphic section of $\rho_{\Delta(\tau_v,i;1),s}$, such that $f_{\Delta(\tau_v,i;1),s_0}=f_0$. Then
$$
\Lambda_{v,s_0}(f_0,\phi_v)=\Lambda(f_{\Delta(\tau_v,i;1),s},\phi_v)_{\large |_{s=s_0}}.
$$
This is well defined, since $\Lambda(f_{\Delta(\tau,i;1)_v,s},\phi_v)_{\large |_{s=s_0}}$ depends only on $f_{\Delta(\tau,i;1)_v,s_0}$. Let $\xi_0$ be an element in the dual of $\rho_{\Delta(\tau_v,i)\gamma_{\psi_v},s_0}$, realized as $\rho_{\hat{\Delta(\tau_v,i)}\gamma^{-1}_{\psi_v},-s_0}$. Assume that it is zero on the image of $\Lambda_{v,s_0}$. Then we need to prove that $\xi_0=0$.
Our assumption is that for all $f_0$ and $\phi_v$,
\begin{equation}\label{(16)}
<\Lambda_{v,s_0}(f_0,\phi_v),\xi_0>=0.
\end{equation}
Let us take in \eqref{(16)}, $f_0=\rho(\hat{w}_{n(i+1)})f'_0$, where $f'_0$ is supported in the open cell, modulo $Q_{ni,n}(F_v)$ from the left, in $\Sp_{2n(i+1)}(F_v)$, and assume that this support is compact. We claim that the support of the function
\begin{equation}\label{(17)}
z\mapsto \Lambda_{v,s_0}(f_0,\phi_v)(J'_{2ni}u_{ni}(z))
\end{equation}
is compact. Here, $S_{ni}(F_v)$ is the space of $ni\times ni$ matrices $z$ over $F_v$, such that $w_{ni}z$ is symmetric. Also,
$J'_{2ni}=\begin{pmatrix}&I_{ni}\\-I_{ni}&\end{pmatrix}$. By \eqref{(6)},
we have\\
\\
$\Lambda_{v,s_0}(f_0,\phi_v)(J'_{2ni}u_{ni}(z))=$
\begin{equation}\label{(18)}
\int \psi^{-1}_v(y_{n,1})\omega_{\psi_v^{-1}}((J_{2ni}u_{ni}(z),1))\phi_v(x_n)f'_0(J_{2n(i+1)}u_{n(i+1)}(\begin{pmatrix}x&y\\z&x'\end{pmatrix}^{w_{n(i+1)}}))dxdy.
\end{equation}
The integration is over $x\in M_{n\times ni}(F_v)$ and $y\in S_n(F_v)$.
By our assumption on $f'_0$, the r.h.s. of \eqref{(18)} is compactly supported in $z$. For such $f'_0$, we can write the l.h.s. of \eqref{(16)} as an absolutely convrgent integral along the open cell, modulo $Q^{(2)}_{ni}(F_v)$ from the left, in $\Sp^{(2)}_{2ni}(F_v)$, and we get
\begin{equation}\label{(19)}
\int_{S_{ni}(F_v)}<\Lambda_{v,s_0}(f_0,\phi_v)(J'_{2ni}u_{ni}(z)),\xi_0(J'_{2ni}u_{ni}(z))>dz=0.
\end{equation}
The inner pairing $<\ ,\ >$ in \eqref{(19)} is the bilinear invariant pairing on $\Delta(\tau_v,i)\times \widehat{\Delta(\tau_v,i)}$. By \eqref{(18)}, we can rewrite \eqref{(19)} as
\begin{multline}\label{(20)}
\int \psi^{-1}_v(y_{n,1})\omega_{\psi_v^{-1}}((J_{2ni}u_{ni}(z),1))\phi_v(x_n)\\
<f'_0(J_{2n(i+1)}u_{n(i+1)}(\begin{pmatrix}x&y\\z&x'\end{pmatrix}^{w_{n(i+1)}})),\xi_0(J'_{2ni}u_{ni}(z))>dxdydz=0.
\end{multline}
The integration is over $z\in S_{ni}(F_v)$, $x\in M_{n\times ni}(F_v)$ and $y\in S_n(F_v)$.
This is valid for all $f'_0$ as above and all $\phi_v$. We conclude that $\xi_0$ must vanish on the open cell, modulo $Q^{(2)}_{ni}(F_v)$ from the left, in $\Sp^{(2)}_{2ni}(F_v)$, and hence $\xi_0=0$. The proposition follows.
	
\end{proof}

\section{Top orbits for $\mathcal{E}^{\Sp_{2mn}}_{\Delta(\tau,m),\wedge^2,[\frac{m+1}{2}]}$ and $\mathcal{E}^{\Sp^{(2)}_{2mn}}_{\Delta(\tau,m)\gamma_\psi,\vee^2,[\frac{m+1}{2}]}$}

We will need to prove the following initial cases in order to prove Theorem \ref{thm 1.2} by induction. In this section, we consider the following residual representations of $\Sp_{2mn}^{(\epsilon)}(\BA)$:\\
1. When $L(\tau,\wedge^2,s)$ has a pole at $s=1$, and $L(\tau,\frac{1}{2})\neq 0$, we know from \cite{JLZ13} (see the remark after Prop. \ref{prop (0)}) that
$E(f_{\Delta(\tau,m),s})$ (on $\Sp_{2mn}(\BA)$) has a simple pole at $s=\frac{m}{2}$, as the section varies. We will consider the residual representation $\mathcal{E}^{\Sp_{2mn}}_{\Delta(\tau,m),\wedge^2,[\frac{m+1}{2}]}$ generated by the residues $Res_{s=\frac{m}{2}}E^{\Sp_{2nm}}(f_{\Delta(\tau,m),s})$. Note that in this case, $n$ must be even and $\omega_\tau=1$.\\
2. When $L(\tau,\vee^2,s)$ has a pole at $s=1$, we know from Prop. \ref{prop (0)} that\\
$E^{\Sp^{(2)}_{2mn}}(f_{\Delta(\tau,m)\gamma_\psi,s})$ has a simple pole at $s=\frac{m}{2}$, as the section varies. We will consider the residual representation $\mathcal{E}^{\Sp^{(2)}_{2mn}}_{\Delta(\tau,m)\gamma_\psi,\vee^2,[\frac{m+1}{2}]}$ generated by the residues\\
 $Res_{s=\frac{m}{2}}E^{\Sp^{(2)}_{2mn}}(f_{\Delta(\tau,m),s})$.

\begin{thm}\label{thm 4.1}
	Assume that $L(\tau,\wedge^2,s)$ has a pole at $s=1$, and $L(\tau,\frac{1}{2})\neq 0$. Then
	$$
	\mathcal{O}(\mathcal{E}^{\Sp_{2mn}}_{\Delta(\tau,m),\wedge^2,[\frac{m+1}{2}]})=(n^{2m}).
	$$
	\end{thm}
The proof of this theorem was sketched by Ginzburg in \cite{G03}, Prop. 3.2. See \cite{L13}, Theorem 1.2, for a detailed proof.

\begin{thm}\label{thm 4.2}
	Assume that $L(\tau,\vee^2,s)$ has a pole at $s=1$.
	\begin{enumerate}
	\item If $\omega_\tau=1$, then
	$$
		\mathcal{O}(\mathcal{E}^{\Sp^{(2)}_{2mn}}_{\Delta(\tau,m)\gamma_\psi,\vee^2,[\frac{m+1}{2}]})=
		\begin{cases}
		(n^{2m}),\ \ \ \ \ n\ even\\
	 ((n+1)^m,(n-1)^m),\ \ n\ odd
	 \end{cases}
	$$
	\item If $\omega_\tau\neq 1$ (recall that $\omega^2_\tau=1$) then
	$$
	\mathcal{O}(\mathcal{E}^{\Sp^{(2)}_{2mn}}_{\Delta(\tau,m)\gamma_\psi,\vee^2,[\frac{m+1}{2}]})=
	\begin{cases}
((n+2)^m, (n-2)^m),\ \ n\ even\\
	((n+1)^m,(n-1)^m),\ \ n\ odd
	\end{cases}
	$$
	\end{enumerate}
\end{thm}

{\bf Proof of Theorem \ref{thm 4.2}, Part 1:}\\
\vspace{0.1cm}

Assume that $\omega_\tau=1$. In this case, the proof of the theorem when $n$ is even is the same as that of Theorem \ref{thm 4.1}. We prove the theorem when $n=2n'-1$ is odd. We first note that $\mathcal{O}(\mathcal{E}^{\Sp^{(2)}_{2mn}}_{\Delta(\tau,m)\gamma_\psi,\vee^2,[\frac{m+1}{2}]})$ is bounded by $((n+1)^m,(n-1)^m)$. This is in \cite{GS20}, Prop. 3.2, where it is stated for $m$ even, but the same proof works for $m$ odd as well. The main work of the proof of Theorem \ref{thm 4.2} in this case is to show that $\mathcal{E}_{\Delta(\tau,m)\gamma_\psi,\vee^2,[\frac{m+1}{2}]}$ admits a nontrivial Fourier coefficient, corresponding to the partition $((n+1)^m,(n-1)^m)$. The corresponding Fourier coefficient (see \cite{GRS03}, Sec. 2) is with respect to the unipotent radical $U_{m,(2m)^{n'-1}}=U^{\Sp_{2mn}}_{m,(2m)^{n'-1}}$ and the character $\psi_{U_{m,(2m)^{n'-1}}}$, as follows. Write an element of $U_{m,(2m)^{n'-1}}(\BA)$ as
\begin{equation}\label{4.1}
u=\begin{pmatrix} I_m&x_{1,2}&\ast&\cdots&\ast&\ast&\ast&\cdots&\ast&\ast\\
&I_{2m}&x_{2,3}&&\ast&\ast&\ast&&&\ast\\&&\ddots&&\vdots&\vdots&&&&\vdots\\&&&I_{2m}&x_{n'-1,n'}&\ast&\ast&&&\ast\\&&&&I_{2m}&x_{n',n'+1}&\ast&&&\ast\\&&&&&I_{2m}&x'_{n'-1,n'}&&&\ast\\&&&&&&I_{2m}&&&\ast\\&&&&&&&\ddots&&\vdots\\&&&&&&&&x'_{2,3}&\ast\\&&&&&&&&I_{2m}&x'_{1,2}\\&&&&&&&&&I_m\end{pmatrix}.
\end{equation}
Then
\begin{multline}\label{4.2}
\psi_{U_{m,(2m)^{n'-1}}}(u)=\\
\psi(tr(x_{1,2}\begin{pmatrix}0_{m\times m}\\I_m\end{pmatrix})+tr(x_{2,3})+tr(x_{3,4})+\cdots tr(x_{n'-1,n'})+\frac{1}{2}tr(x_{n',n'+1})).
\end{multline}
Note that $x_{n',n'+1}$ has the form
$$
x_{n',n'+1}=\begin{pmatrix}a&b\\c&w_m{}^taw_m\end{pmatrix},
$$
where $a,b,c$ are $m\times m$ matrices, so that $\frac{1}{2}tr(x_{n',n'+1})=tr(a)$.

Our goal is to show that
\begin{equation}\label{4.3}
\mathcal{F}_{\psi_{U_{m,(2m)^{n'-1}}}}(\xi)=\int_{U_{m,(2m)^{n'-1}}(F)\backslash U_{m,(2m)^{n'-1}}(\BA)}\xi(u)\psi^{-1}_{U_{m,(2m)^{n'-1}}}(u)du\not\equiv 0,
\end{equation}
as $\xi$ varies in $\mathcal{E}_{\Delta(\tau,m)\gamma_\psi,\vee^2,[\frac{m+1}{2}]}$. Recall that we identify $U_{m,(2m)^{n'-1}}(\BA)$ as a subgroup of $\Sp_{2mn}^{(2)}(\BA)$.
Consider the following Weyl element $w_0\in \Sp_{2nm}(F)$. Write $w_0$ as a $2n\times 2n$ matrix of $m\times m$ blocks. Then, for $1\leq i\leq n$, the $i$-th block row of $w_0$ has the form
$$
\begin{pmatrix}0_{m\times m}&0_{m\times m}&\cdots&0_{m\times m}&I_m&0_{m\times m}&\cdots&0_{m\times m}\end{pmatrix},
$$
where $I_m$ appears at the $2i-1$ position. This determines $w_0$. For example, for $n=3$
$$
w_0=\begin{pmatrix}I_m&0&0&0&0&0\\0&0&I_m&0&0&0\\0&0&0&0&I_m&0\\0&-I_m&0&0&0&0&\\0&0&0&I_m&0&0\\0&0&0&0&0&I_m\end{pmatrix}.
$$
We may conjugate inside the integral \eqref{4.3} by $w_0$,
\begin{multline}\label{4.4}
\mathcal{F}_{\psi_{U_{m,(2m)^{n'-1}}}}(\xi)=\int_{U_{m,(2m)^{n'-1}}(F)\backslash U_{m,(2m)^{n'-1}}(\BA)}\xi(w_0u)\psi^{-1}_{U_{m,(2m)^{n'-1}}}(u)du=\\
=\int_{V(F)\backslash V(\BA)}\xi(vw_0)\psi^{-1}_V(v)dv,
\end{multline}
where $V=w_0U_{m,(2m)^{n'-1}}w_0^{-1}$, and, for $v\in V(\BA)$,
	$$
	\psi_V(v)=\psi_{U_{m,(2m)^{n'-1}}}(w_0^{-1}vw_0).
	$$
Let us describe the subgroup $V$. Write $v\in V$ as
\begin{equation}\label{4.5}
v=\begin{pmatrix}U&X\\Y&U'\end{pmatrix}.
\end{equation}
Then $U$ has the form
\begin{equation}\label{4.6}
U=\begin{pmatrix}I_m&u_{1,2}&\ast&\cdots&\ast\\&I_m&u_{2,3}&&\ast\\&&\ddots&&\vdots\\&&&I_m&u_{n-1,n}\\&&&&I_m\end{pmatrix};
\end{equation}
the matrix $U'$ has a similar form to \eqref{4.6}. Next, write the matrices $X$, $Y$, each, as an $n\times n$ matrix of $m\times m$ blocks. Then $X$ and $Y$ have upper triangular shapes. The matrix $Y$ is also such that its main block diagonal and the one above it consist of zero matrices. Thus, $X$, $Y$ have the forms
\begin{multline}\label{4.7}
X=\begin{pmatrix}a_{1,1}&a_{1,2}&\ast&\cdots&a_{1,n}\\&a_{2,2}&a_{2,3}&&a_{2,n}\\&&\ddots&&\vdots\\&&&a_{n-1,n-1}&a_{n-1,n}\\&&&&a_{n,n}\end{pmatrix};\\
Y=\begin{pmatrix}0_{m\times m}&0_{m\times m}&b_{1,3}&b_{1,4}&\cdots&b_{1,n}\\&0_{m\times m}&0_{m\times m}&b_{2,4}&&b_{2,n}\\&&\ddots&&&\vdots\\&&&0_{m\times m}&0_{m\times m}&b_{n-2,n}\\&&&&0_{m\times m}&0_{m\times m}\\&&&&&0_{m\times m}\end{pmatrix}.
\end{multline}	
For $v\in V(\BA)$ of the form \eqref{4.5}-\eqref{4.7},
\begin{equation}\label{4.8}
\psi_V(v)=\psi(tr(u_{1,2}+u_{2,3}+\cdots u_{n',n'+1})-tr(u_{n'+1,n'+2}+u_{n'+2,n'+3}+\cdots+u_{n-1,n})).
\end{equation}	
Note that we have here a special case of (5.8) in \cite{GS20}, where in (5.9) in loc. cit. we should replace $r$ by $m$, $n$ by $n'$, and we should also ignore in (5.8) there the second block row and block column.
Now we apply a chain of roots exchange, exactly as we did in the second part of the proof of Prop. 6.1 in \cite{GS20}, (where we exchanged there $\RY_{3,1}^{i-1,j}$ with $\RX_{1,3}^{j-1,i-1}$, for $1<i\leq j-2$, and $\RY_{1,3}^{j-1,j+1}$ with $\RX_{1,3}^{j,j-1}$). Thus, we exchange the $(1,3)$ block of $Y$ with the $(2,1)$ block of $X$, then the $(2,4)$, $(1,4)$ blocks of $Y$ with the $(3,2)$, $(3,1)$ blocks of $X$, in this order, and we proceed and exchange, for $1\leq j\leq n'$, the $(j-2,j)$, $(j-3,j)$,..., $(1,j)$ blocks of $Y$ with the $(j-1,j-2)$, $(j-1,j-3)$,..., $(j-1,1)$ blocks of $X$, in this order. We get, as in Prop. 6.2 in \cite{GS20}, that the Fourier coefficient \eqref{4.3} (and hence \eqref{4.4}), $\mathcal{F}_{\psi_{U_{m,(2m)^{n'-1}}}}$, is nontrivial on $\mathcal{E}_{\Delta(\tau,m)\gamma_\psi,\vee^2,[\frac{m+1}{2}]}$, if and only if the following Fourier coefficient is nontrivial on $\mathcal{E}_{\Delta(\tau,m)\gamma_\psi,\vee^2,[\frac{m+1}{2}]}$,   	
\begin{equation}\label{4.9}
\mathcal{F}_{\psi_{V^{n',n'}}}(\xi)=\int_{V^{n',n'}(F)\backslash V^{n',n'}(\BA)}\xi(u)\psi^{-1}_{V^{n',n'}}(u)du,
\end{equation}	
where $V^{n',n'}$ is the subgroup of symplectic matrices of the form \eqref{4.5}, with $U$ (and $U'$) of the form \eqref{4.6}, but now $X$, $Y$ are of the following forms,
\begin{equation}\label{4.10}
X=\begin{pmatrix} a_{1,1}&\cdots&a_{1,n'-1}&a_{1,n'}&a_{1,n'+1}&\cdots&a_{1,n}\\\vdots&&&\vdots&&&\vdots\\a_{n'-1,1}&\cdots&a_{n'-1,n'-1}&a_{n'-1,n'}&a_{n'-1,n'+1}&\cdots&\ast\\0&\cdots&0&a_{n',n'}&\ast&&\ast\\&&&0&\ast&&\ast\\&&&\vdots&\vdots&&\vdots\\&&&0&\ast&&\ast\end{pmatrix}
\end{equation}
Recall that $X$ is written as a $n\times n$ matrix of $m\times m$ blocks.
\begin{equation}\label{4.11}
Y=\begin{pmatrix}0&0&\cdots&0&b_{1,n'+1}&b_{1,n'+2}&\cdots&b_{1,n-1}&b_{1,n}\\&0&&0&b_{2,n'+1}&b_{2,n'+2}&&b_{2,n-1}&\ast\\&&\ddots&&&&&\vdots\\&&&0&b_{n'-1,n'+1}&\ast&&\ast&\ast\\&&&&0&0&\cdots&0&0\\&&&&&\ddots&&&\vdots\\&&&&&&&0&0\\&&&&&&&&0\end{pmatrix}.
\end{equation}
The character $\psi_{V^{n',n'}}$ is given by the same formula as in \eqref{4.8}.	

Next, we continue exactly as in the second step, right after Prop. 6.2 in \cite{GS20}, meaning that we exchange the $(n'-1,n'+1)$ block of $Y$ "into" the $(n',n'-1)$ block of $X$. Note that if we let $\RY^{n'-1,n'+1}$ denote the subgroup of elements of the form \eqref{4.5}, with $U=I_{mn}$, $X=0$ and $Y$ as in \eqref{4.11}, with all blocks being zero, except for $b_{n'-1,n'+1}$, then $w_mb_{n'-1,n'+1}$ must be a symmetric matrix. If we similarly define the subgroup $\RX^{n'-1,n'}$, then it isomorphic to the additive group of $m\times m$ matrices (as algebraic groups over $F$). Thus, we exchange $\RY^{n'-1,n'+1}$ with the subgroup of $\RX^{n'-1,n'}$, corresponding to the $m\times m$ matrices $x$, such that $w_mx$ is symmetric. It follows that the Fourier coefficient \eqref{4.9} is nontrivial on $\mathcal{E}_{\Delta(\tau,m)\gamma_\psi,\vee^2,[\frac{m+1}{2}]}$, if and only if the following Fourier coefficient is nontrivial on $\mathcal{E}_{\Delta(\tau,m)\gamma_\psi,\vee^2,[\frac{m+1}{2}]}$,   	
\begin{equation}\label{4.12}
\mathcal{F}_{\psi_{V^{n',n'-1;0}}}(\xi)=\int_{V^{n',n'-1;0}(F)\backslash V^{n',n'-1;0}(\BA)}\xi(u)\psi^{-1}_{V^{n',n'-1;0}}(u)du,
\end{equation}	
where $V^{n',n'-1;0}$ is the subgroup of symplectic matrices of the form \eqref{4.5}, with $U$ (and $U'$) of the form \eqref{4.6}, $Y$ has the form \eqref{4.11} with $b_{n'-1,n'+1}=0$ and $X$ has the form,
\begin{equation}\label{4.13}
X=\begin{pmatrix} a_{1,1}&a_{1,2}&\cdots&a_{1,n'-1}&a_{1,n'}&a_{1,n'+1}&\cdots&a_{1,n-1}&a_{1,n}\\\vdots&&&&\vdots&&&&\vdots\\a_{n'-1,1}&a_{n'-1,2}&\cdots&a_{n'-1,n'-1}&a_{n'-1,n'}&a_{n'-1,n'+1}&\cdots&\ast&\ast\\0&0&\cdots&a_{n',n'-1}&a_{n',n'}&\ast&&\ast&\ast\\&&&&\ast&\ast&&\ast&\ast\\&&&&\vdots&\vdots&&&\vdots\\&&&&0&\ast&&\ast&\ast\end{pmatrix},
\end{equation}	
such that $w_ma_{n',n'-1}$ is symmetric. The character $\psi_{V^{n',n'-1;0}}$ is given by the same formula as in \eqref{4.8}.	Now, we apply the argument of Prop. 7.1 in \cite{GS20} and get that in \eqref{4.12}, we may replace $V^{n',n'-1;0}$ by $V^{n',n'-1}$ which is defined as $V^{n',n'-1;0}$, except that in \eqref{4.13}, we do not restrict $w_ma_{n',n'-1}$ to be symmetric, and let it be arbitrary. The character $\psi_{V^{n',n'-1}}$ is given by the same formula as in \eqref{4.8}. The point is that when we carry out the Fourier expansion as in Prop. 7.1 in \cite{GS20}, we get Fourier coefficients corresponding to partitions of the form $(n+2,...)$.
These must be zero on $\mathcal{E}_{\Delta(\tau,m)\gamma_\psi,\vee^2,[\frac{m+1}{2}]}$, since $\mathcal{O}(\mathcal{E}_{\Delta(\tau,m)\gamma_\psi,\vee^2,[\frac{m+1}{2}]})$ is bounded by $((n+1)^m,(n-1)^m)$. Thus, the Fourier coefficient \eqref{4.12} is nontrivial on $\mathcal{E}_{\Delta(\tau,m)\gamma_\psi,\vee^2,[\frac{m+1}{2}]}$, if and only if the following Fourier coefficient is nontrivial on $\mathcal{E}_{\Delta(\tau,m)\gamma_\psi,\vee^2,[\frac{m+1}{2}]}$,   	
\begin{equation}\label{4.14}
\mathcal{F}_{\psi_{V^{n',n'-1}}}(\xi)=\int_{V^{n',n'-1}(F)\backslash V^{n',n'-1}(\BA)}\xi(u)\psi^{-1}_{V^{n',n'-1}}(u)du.
\end{equation}
Now, we exchange roots, as we did right after Prop. 7.1 in \cite{GS20}, namely exchange the $(n'-2,n'+1)$, $(n'-3,n'+1)$,...,$(1,n'+1)$ blocks in $Y$ with the blocks in $X$ in positions $(n',n'-2)$, $(n',n'-3)$,...,$(n',1)$. We get that the Fourier coefficient \eqref{4.14} is nontrivial on $\mathcal{E}_{\Delta(\tau,m)\gamma_\psi,\vee^2,[\frac{m+1}{2}]}$, if and only if the following Fourier coefficient is nontrivial on $\mathcal{E}_{\Delta(\tau,m)\gamma_\psi,\vee^2,[\frac{m+1}{2}]}$,   	
\begin{equation}\label{4.15}
\mathcal{F}_{\psi_{V^{n',1}}}(\xi)=\int_{V^{n',1}(F)\backslash V^{n',1}(\BA)}\xi(u)\psi^{-1}_{V^{n',1}}(u)du,
\end{equation} 	
where $V^{n',1}$ is as $V^{n',n'-1}$ except that the block $Y$, written as in \eqref{4.11}, is such that all blocks in column $n'+1$ (and row $n'-1$) are zero, and in $X$, written as in \eqref{4.10}, we fill in the zero blocks in row $n'$ (and column $n'$).	The character $\psi_{V^{n',1}}$ is given by the same formula as in \eqref{4.8}. We continue like this, step by step, repeating the arguments of Prop. 7.4, 7.1 in \cite{GS20}, and then we exchange roots. We write Fourier expansions, showing at each step that only the trivial character contributes to the Fourier expansion. At each such step, we use the fact that $\mathcal{O}(\mathcal{E}_{\Delta(\tau,m)\gamma_\psi,\vee^2,[\frac{m+1}{2}]})$ is bounded by $((n+1)^m,(n-1)^m)$. In the end we get that the Fourier coefficient \eqref{4.4} is nontrivial on $\mathcal{E}_{\Delta(\tau,m)\gamma_\psi,\vee^2,[\frac{m+1}{2}]}$, if and only if the following Fourier coefficient is nontrivial on $\mathcal{E}_{\Delta(\tau,m)\gamma_\psi,\vee^2,[\frac{m+1}{2}]}$,   	
\begin{equation}\label{4.16}
\mathcal{F}_{\psi_{U_{m^n}}}(\xi)=\int_{U_{m^n}(F)\backslash U_{m^n}(\BA)}\xi(u)\psi^{-1}_{U_{m^n}}(u)du,
\end{equation} 	
where, for
$$
v=\begin{pmatrix}U&X\\&U^*\end{pmatrix}\in U_{m^n}(\BA),
$$
with $U$ as in \eqref{4.6},	$\psi_{U_{m^n}}(v)$ is given by \eqref{4.8}. Thus, $\mathcal{F}_{\psi_{U_{m^n}}}$ is the composition of the constant term of $\xi$, $\xi^{U_{mn}}$, along $U_{mn}$, and the Fourier coefficient along the unipotent radical $V_{m^n}$ (inside $\GL_{mn}$), and the character given by \eqref{4.8}, which we denote now by $\psi_{V_{m^n}}$. From the proof of Prop. \ref{prop (0)}, \eqref{4.16} is nontrivial on\\
$\mathcal{E}_{\Delta(\tau,m)\gamma_\psi,\vee^2,[\frac{m+1}{2}]}$. This proves the first part of the theorem.\\
\vspace{0.1cm}

{\bf Proof of Theorem \ref{thm 4.2}, Part 2:}\\
\vspace{0.1cm}

Assume now, that the quadratic character $\omega_\tau$ is nontrivial. Let $v$ be a finite place, where $\tau_v$ is unramified, and its central character is the unique, unramified, nontrivial character $\lambda_v$. Assume that $n=2n'$ is even. Since $\tau_v$ is self dual, we can write $\tau_v$ as a parabolic induction from an unramified character of the standard Borel subgroup, as follows,
\begin{equation}\label{4.19.1}
\tau_v=\chi_1\times\cdots\times \chi_{n'-1}\times 1\times\lambda_v\times \chi^{-1}_{n'-1}\times\cdots\times\chi^{-1}_1,
\end{equation}
where $\chi_i$ are unramified characters of $F_v^*$. Then the unramified constituent of the factor at $v$ of 	
 $\mathcal{E}_{\Delta(\tau,m)\gamma_\psi,\vee^2,[\frac{m+1}{2}]}$  is the unramified constituent of the following parabolic induction,
\begin{equation}\label{4.19.2}
\Ind^{\Sp^{(2)}_{2mn}(F_v)}_{Q^{(2)}_{(2m)^{n'-1},m^2}(F_v)}[\otimes_{i=1}^{n'-1}(\chi_i\circ det_{\GL_{2m}}) \otimes (\lambda_v\circ det_{\GL_m})|det_{\GL_m}|^{\frac{m}{2}}\otimes |det_{\GL_m}|^{\frac{m}{2}}]\gamma_{\psi_v}.
\end{equation}
Recall that the pole in question here is at $s=\frac{m}{2}$. As in the proof of Prop. 3.2 in \cite{GS20}, we conclude, from \eqref{4.19.2}, that all symplectic partitions of $2mn$ corresponding to nontrivial Fourier coefficients on  $\mathcal{E}^{\Sp^{(2)}_{2mn}}_{\Delta(\tau,m)\gamma_\psi,\vee^2,[\frac{m+1}{2}]}$ are bounded by the induced nilpotent orbit corresponding to \eqref{4.19.2}, and this  corresponds to the partition $((n+2)^m, (n-2)^m))$ of $2nm$. See
\cite{CM93}, Chapter 7. When $n=2n'+1$ is odd, the analog of \eqref{4.19.1} is, with similar notation,
\begin{equation}\label{4.19.3}
\tau_v=\chi_1\times\cdots\times \chi_{n'}\times\lambda_v\times \chi^{-1}_{n'}\times\cdots\times\chi^{-1}_1,
\end{equation}
and the analog of \eqref{4.19.2} is
\begin{equation}\label{4.19.4}
\Ind^{\Sp^{(2)}_{2mn}(F_v)}_{Q^{(2)}_{(2m)^{n'},m}(F_v)}[\otimes_{i=1}^{n'}(\chi_i\circ det_{\GL_{2m}}) \otimes (\lambda_v\circ det_{\GL_m})|det_{\GL_m}|^{\frac{m}{2}}]\gamma_{\psi_v}.
\end{equation}
In this case, as in the case where $n$ is odd and the central character is trivial, the associated partition for \eqref{4.19.4} is $((n+1)^m,(n-1)^m)$, as in the first part of the theorem. What we proved in the first part works exactly the same in this case. This proves the second part of the theorem when $n$ is odd. Thus, assume that $n=2n'$ is even (and $\omega_\tau\neq 1$).
As in the first part, we need to show that $\mathcal{E}_{\Delta(\tau,m)\gamma_\psi,\vee^2,[\frac{m+1}{2}]}$ admits a (nontrivial) Fourier coefficient, corresponding to the partition $((n+2)^m,(n-2)^m)$. The corresponding Fourier coefficient (see \cite{GRS03}, Sec. 2) is with respect to the unipotent radical $U_{m^2,(2m)^{n'-1}}=U^{\Sp_{2mn}}_{m^2,(2m)^{n'-1}}$ and the character $\psi_{U_{m^2,(2m)^{n'-1}}}$ as follows. Write an element of $U_{m^2,(2m)^{n'-1}}(\BA)$ as
\begin{equation}\label{4.20}
e=\begin{pmatrix}I_m&x&z\\&u&x'\\&&I_m\end{pmatrix}\in \Sp_{2mn}(\BA),
\end{equation}
where $u$ has the form \eqref{4.1}. We keep using the notation in \eqref{4.1}, only that now, with $n'=\frac{n}{2}$, $u$ lies in $\Sp_{2m(n-1)}(\BA)$. Write in \eqref{4.20},
$$
x=(x_{0,1}, x_{0,2},...,x_{0,n}),
$$
where $x_{0,1}, x_{0,n}\in M_{m\times m}(\BA)$,  $x_{0,2},...,x_{0,n-1}\in M_{m\times 2m}(\BA)$.
The character \\
$\psi_{U_{m^2,(2m)^{n'-1}}}$ is given by
\begin{multline}\label{4.21}
\psi_{U_{m^2,(2m)^{n'-1}}}(u)=\\
\psi(tr(x_{0,1})+tr(x_{1,2}\begin{pmatrix}I_m\\0_{m\times m}\end{pmatrix})+tr(x_{2,3})+tr(x_{3,4})+\cdots+tr(x_{n'-1,n'})+\frac{1}{2}tr(x_{n',n'+1})).
\end{multline}
As in \eqref{4.3}, our goal is to show that
\begin{equation}\label{4.22}
\mathcal{F}_{\psi_{U_{m^2,(2m)^{n'-1}}}}(\xi)=\int_{U_{m^2,(2m)^{n'-1}}(F)\backslash U_{m^2,(2m)^{n'-1}}(\BA)}\xi(u)\psi^{-1}_{U_{m^2,(2m)^{n'-1}}}(u)du\not\equiv 0,
\end{equation}
as $\xi$ varies in $\mathcal{E}^{\Sp^{(2)}_{2nm}}_{\Delta(\tau,m)\gamma_\psi,\vee^2,[\frac{m+1}{2}]}$. The proof is very similar to the first part, and so we will just sketch it.
Consider the following Weyl element $\tilde{w}_0\in \Sp_{2nm}(F)$. Write $\tilde{w}_0$ as a $2n\times 2n$ matrix of $m\times m$ blocks. Then the first three block rows form the matrix
$$
\begin{pmatrix}I_{3m}&0&\cdots&0\end{pmatrix}.
$$
For $4\leq i\leq n$, the $i$-th block row of $\tilde{w}_0$ has the form
$$
\begin{pmatrix}0_{m\times m}&0_{m\times m}&\cdots&0_{m\times m}&I_m&0_{m\times m}&\cdots&0_{m\times m}\end{pmatrix},
$$
where $I_m$ appears at the $2i-3$ position. This determines $\tilde{w}_0$. For example, for $n=6$
$$
\tilde{w}_0=\begin{pmatrix}I_{3m}&0&0&0&0&0&0&0\\0&0&I_m&0&0&0&0&0\\0&0&0&0&I_m&0&0&0\\ 0&0&0&0&0&0&I_m&0\\0&-I_m&0&0&0&0&0&0\\0&0&0&-I_m&0&0&0&0\\0&0&0&0&0&I_m&0&0\\0&0&0&0&0&0&0    &I_{3m}\end{pmatrix}.
$$
We conjugate inside the integral \eqref{4.22} by $\tilde{w}_0$,
\begin{multline}\label{4.23}
\mathcal{F}_{\psi_{U_{m^2,(2m)^{n'-1}}}}(\xi)=\int_{U_{m^2,(2m)^{n'-1}}(F)\backslash U_{m^2,(2m)^{n'-1}}(\BA)}\xi(\tilde{w}_0u)\psi^{-1}_{U_{m^2,(2m)^{n'-1}}}(u)du=\\
=\int_{V(F)\backslash V(\BA)}\xi(vw_0)\psi^{-1}_V(v)dv,
\end{multline}
where $V=\tilde{w}_0U_{m^2,(2m)^{n'-1}}\tilde{w}_0^{-1}$, and, for $v\in V(\BA)$,
$$
\psi_V(v)=\psi_{U_{m^2,(2m)^{n'-1}}}(\tilde{w}_0^{-1}v\tilde{w}_0).
$$
The subgroup $V$ has a similar description to \eqref{4.5} - \eqref{4.7}. In the notation of \eqref{4.5}, $U$ has the form \eqref{4.6}. Write the matrices $X$, $Y$, each, as $n\times n$ matrices of $m\times m$ blocks. They have the forms
\begin{multline}\label{4.24}
X=\begin{pmatrix}a_{1,1}&a_{1,2}&a_{1,3}&\cdots&a_{1,n-1}&a_{1,n}\\a_{2,1}&a_{2,2}&a_{2,3}&&a_{2,n-1}&a_{2,n}\\&a_{3,2}&a_{3,3}&&a_{3,n-1}&a_{3,n}\\&&\ddots&\ddots&\vdots&\vdots\\&&&a_{n-1,n-2}&a_{n-1,n-1}&a_{n-1,n}\\&&&&a_{n-1,n}&a_{n,n}\end{pmatrix};\\
Y=\begin{pmatrix}0_{m\times m}&0_{m\times m}&0_{m\times m}&b_{1,4}&b_{1,5}&\cdots&b_{1,n}\\&0_{m\times m}&0_{m\times m}&0_{m\times m}&b_{2,5}&&b_{2,n}\\&&\ddots&&&\vdots\\&&&0_{m\times m}&0_{m\times m}&0_{m\times m}&b_{n-3,n}\\&&&&0_{m\times m}&0_{m\times m}&0_{m\times m}\\&&&&&0_{m\times m}&0_{m\times m}\\&&&&&&0_{m\times m}\end{pmatrix}.
\end{multline}
The matrix $Y$ has an upper triangular shape, with its main block diagonal and two diagonals above it being zero. The matrix $X$ is such that its lower $n-2$ block diagonals are zero. Finally, for $v\in V(\BA)$ of the form \eqref{4.5}, \eqref{4.6}, \eqref{4.24},
\begin{equation}\label{4.25}
\psi_V(v)=\psi(tr(u_{1,2}+u_{2,3}+\cdots u_{n'+1,n'+2})-tr(u_{n'+2,n'+2}+u_{n'+3,n'+4}+\cdots+u_{n-1,n})).
\end{equation}
Now, we exchange roots "from $Y$ into $X$" as we did in the first part. We exchange the block $b_{1,4}$ in $Y$ (notation of \eqref{4.24}) to "fill" the block in $X$ in position $(3,1)$, then the blocks $b_{2,5}$, $b_{1,5}$ in $Y$ to fill in the blocks in $X$ in positions $(4,1)$, $(4,2$ (in this order), and so on, exactly as in the first part. The analogs of \eqref{4.10}, \eqref{4.11} are as follows.
\begin{equation}\label{4.26}
X=\begin{pmatrix} a_{1,1}&\cdots&a_{1,n'-1}&a_{1,n'}&a_{1,n'+1}&\cdots&a_{1,n}\\\vdots&&&\vdots&&&\vdots\\a_{n',1}&\cdots&a_{n',n'-1}&a_{n',n'}&a_{n',n'+1}&\cdots&\ast\\0&\cdots&0&a_{n'+1,n'}&\ast&&\ast\\&&&0&\ast&&\ast\\&&&\vdots&\vdots&&\vdots\\&&&0&\ast&&\ast\end{pmatrix}
\end{equation}

\begin{equation}\label{4.27}
Y=\begin{pmatrix}0&0&\cdots&0&b_{1,n'+2}&b_{1,n'+3}&\cdots&b_{1,n-1}&b_{1,n}\\&0&&0&b_{2,n'+2}&b_{2,n'+3}&&b_{2,n-1}&\ast\\&&\ddots&&&&&\vdots\\&&&0&b_{n'-1,n'+2}&\ast&&\ast&\ast\\&&&&0&0&\cdots&0&0\\&&&&&\ddots&&&\vdots\\&&&&&&&0&0\\&&&&&&&&0\end{pmatrix}.
\end{equation}
We continue as in the first part, using repeatedly the arguments in the proofs of Prop. 6.2, Prop. 7.1, Prop. 7.4 in \cite{GS20}, and the fact that $\mathcal{O}(\mathcal{E}_{\Delta(\tau,m)\gamma_\psi,\vee^2,[\frac{m+1}{2}]})$ is bounded by $((n+2)^m,(n-2)^m)$. In the end we get that the Fourier coefficient \eqref{4.23} is nontrivial on $\mathcal{E}_{\Delta(\tau,m)\gamma_\psi,\vee^2,[\frac{m+1}{2}]}$, if and only if the Fourier coefficient \eqref{4.16} is nontrivial on $\mathcal{E}_{\Delta(\tau,m)\gamma_\psi,\vee^2,[\frac{m+1}{2}]}$, and this we proved in the end of the first part. This completes the proof of the second part of the theorem.

\section{Proof of Theorem \ref{thm 1.2}}

\noindent	
I. Assume that $L(\tau,\wedge^2,s)$ has a pole at $s=1$ and $L(\tau,\frac{1}{2})\neq 0$.\\

Let us write the statement of the theorem in detail for this case. First, the set of points \eqref{1.2} $e^{(2)}_{k,m}(\wedge^2)$ is the set of poles of $E(f_{\Delta(\tau,m)\gamma_\psi,s})$, as the section varies, in $Re(s)\geq 0$, and they are all simple poles. (Recall that we know from \cite{JLZ13} the analogous fact in the linear case ($\epsilon=1$), namely that the set of points \eqref{1.1} $e_{k,m}(\wedge^2)$ is the set of poles of $E(f_{\Delta(\tau,m),s}$, as the section varies, in $Re(s)\geq 0$, and are all simple.) Next,
\begin{enumerate}
	\item For $m$ even, $1\leq k\leq \frac{m}{2}$ ($e_{k,m}^{\Sp_{2mn}}(\wedge^2)=k$)
$$
\mathcal{O}(\mathcal{E}^{\Sp_{2mn}}_{\Delta(\tau,m),\wedge^2,k})=((2n)^{m-2k},n^{4k}).
$$
\item For $m$ odd, $1\leq k\leq \frac{m+1}{2}$, ($e_{k,m}^{\Sp_{2mn}}(\wedge^2)=k-\frac{1}{2}$)
$$
\mathcal{O}(\mathcal{E}^{\Sp_{2mn}}_{\Delta(\tau,m),\wedge^2,k})=((2n)^{m-2k+1},n^{4k-2}).
$$
\item For $m$ even, $1\leq k\leq \frac{m}{2}$, ($e_{k,m}^{\Sp^{(2)}_{2mn}}(\wedge^2)=k-\frac{1}{2}$)
$$
\mathcal{O}(\mathcal{E}^{\Sp^{(2)}_{2mn}}_{\Delta(\tau,m)\gamma_\psi,\wedge^2,k})=((2n)^{m-2k+1},n^{4k-2}).
$$
\item For $m$ odd, $1\leq k\leq \frac{m-1}{2}$, ($e_{k,m}^{\Sp^{(2)}_{2mn}}(\wedge^2)=k$)
$$
\mathcal{O}(\mathcal{E}^{\Sp^{(2)}_{2mn}}_{\Delta(\tau,m)\gamma_\psi,\wedge^2,k})=((2n)^{m-2k},n^{4k}).
$$
\end{enumerate}

Assume that $L(\tau,\wedge^2,s)$ has a pole at $s=1$, and $L(\tau,\frac{1}{2})\neq 0$. We first address the poles of the Eisenstein series $E^{\Sp^{(2)}_{4ni}}(f_{\Delta(\tau,2i)\gamma_\psi,s})$. By Theorem \ref{thm 7.1}, for $h\in \Sp_{2n(2i-1)}(\BA)$,	
\begin{equation}\label{5.1'}
\mathcal{D}^\phi_{\psi,n(2i-1)}(E^{\Sp^{(2)}_{4ni}}(f_{\Delta(\tau,2i)\gamma_\psi,s}))(h)=E^{\Sp_{2n(2i-1)}}(\Lambda(f_{\Delta(\tau,2i)\gamma_\psi,s},\phi))(h).
\end{equation}
Let $1\leq j\leq i$. Note that all the normalizing factors involved in Theorem \ref{thm 7.1} are holomorphic and nonzero at $s=\frac{2j-1}{2}$. From Cor. \ref{cor (1')}, it follows that $\Lambda(f_{\Delta(\tau,2i)\gamma_\psi,s},\phi)$ is holomorphic at $s=\frac{2j-1}{2}$. From Prop. \ref{prop (2)}, it follows that the map\\
 $\Lambda(f_{\Delta(\tau,2i)\gamma_\psi,\frac{2j-1}{2}},\phi)$ on the space of $\rho_{\Delta(\tau,2i)\gamma_\psi,\frac{2j-1}{2}}\times \mathcal{S}(\BA^{(2i-1)n})$ has image which corresponds to
$$
(\otimes_{v\in S_\infty}W(\rho_{\Delta(\tau_v,2i-1),\frac{2j-1}{2}}))\otimes(\otimes'_{v<\infty}V(\rho_{\Delta(\tau_v,2i-1),\frac{2j-1}{2}})),
$$
where $V(\rho_{\Delta(\tau_v,2i-1),\frac{2j-1}{2}})$ denotes the space of $\rho_{\Delta(\tau_v,2i-1),\frac{2j-1}{2}}$, and, for $v\in S_\infty$, $W(\rho_{\Delta(\tau_v,2i-1),\frac{2j-1}{2}})$ is a dense subspace of $V(\rho_{\Delta(\tau_v,2i-1),\frac{2j-1}{2}})$. Taking residues of Eisenstein series at a given point is a continuous map. We conclude from \eqref{5.1'} that $E^{\Sp_{2n(2i-1)}}(\Lambda(f_{\Delta(\tau,2i)\gamma_\psi,s},\phi))$ has a pole at $s=\frac{2j-1}{2}$, for all $1\leq j\leq i$. We used here \cite{JLZ13}.
Thus, the r.h.s. of \eqref{5.1'} has a pole at each point $s=\frac{2j-1}{2}$, $1\leq j\leq i$, and hence, so does 	
$E^{\Sp^{(2)}_{4ni}}(f_{\Delta(\tau,2)\gamma_\psi,s})$ (as the section varies).

By Theorem \ref{thm 7.1}, for all $h\in \Sp^{(2)}_{4ni}(\BA)$,
\begin{equation}\label{5.1''}
\mathcal{D}^\phi_{\psi,2ni}(E^{\Sp_{2n(2i+1)}}(f_{\Delta(\tau,2i+1),s}))(h)=E^{\Sp^{(2)}_{4ni}}(\Lambda(f_{\Delta(\tau,2i+1),s},\phi))(h).
\end{equation}
From what we just proved, we conclude that the r.h.s. of \eqref{5.1''} has poles at all points $\frac{2j-1}{2}$, $1\leq j\leq i$. Again, this follows by the argument above, using Cor. \ref{cor (1')} and Prop. \ref{prop (2)}. Since each point $s=\frac{2j-1}{2}$, $1\leq j\leq i$, is a simple pole of $E^{\Sp_{2n(2i+1)}}(f_{\Delta(\tau,2i+1),s})$, we conclude from \eqref{5.1''}, and the argument above using Cor. \ref{cor (1')} and Prop. \ref{prop (2)}, that the poles at $s=\frac{2j-1}{2}$, $1\leq j\leq i$ of $E^{\Sp^{(2)}_{4ni}}(f_{\Delta(\tau,2)\gamma_\psi,s})$ are all simple (as the section varies). Moreover, let $s_0$, with $Re(s_0)\geq 0$ be a pole of $E^{\Sp^{(2)}_{4ni}}(f'_{\Delta(\tau,2i)\gamma_\psi,s})$. We conclude from \eqref{5.1''} that $s_0$ is a pole of $E^{\Sp_{2n(2i+1)}}(f_{\Delta(\tau,2i+1),s})$, as the section varies. From \cite{JLZ13}, we get that $s_0=\frac{2j-1}{2}$, for some $1\leq j\leq i+1$. We must have $1\leq j\leq i$, since $s=\frac{2i+1}{2}$ is not a pole of $E^{\Sp^{(2)}_{4ni}}(f_{\Delta(\tau,2i)\gamma_\psi,s})$. This follows from the proof of Prop. \ref{prop (0)}.

We will now prove Part I.2 of the theorem, by induction on $m$. The proof then implies Part I.3. Note that the indicated partition, in both cases, $((2n)^{m-2k+1},n^{4k-2})$, bounds $\mathcal{O}(\mathcal{E}^{\Sp^{(\epsilon)}_{2mn}}_{\Delta(\tau,m)\gamma^{(\epsilon)}_\psi,\wedge^2,k})$. This is proved in \cite{GS20}, Prop. 3.1, when $m$ is even. The proof when $m$ is odd is entirely the same.
Thus, we need to show that the residual Eisenstein series admits in each case a nontrivial Fourier coefficient corresponding to the indicated top partition. When an automorphic representation $\pi$ admits a nontrivial Fourier coefficient corresponding to a partition $\underline{p}$, we will also say that $\underline{p}$ supports $\pi$.

We start with $E^{\Sp_{2n}}(f_{\tau,s})$. It has a simple pole at $s=\frac{1}{2}$, and by Theorem \ref{thm 4.1},
\begin{equation}\label{5.1}
\mathcal{O}(\mathcal{E}^{\Sp_{2n}}_{\tau,\wedge^2,1})=(n^2).
\end{equation}
This proves part I.2 of the theorem for $m=1$. Assume by induction that, for $m=2i-1$ odd, and $1\leq k\leq \frac{m+1}{2}=i$,
\begin{equation}\label{5.2}
\mathcal{O}(\mathcal{E}^{\Sp_{2mn}}_{\Delta(\tau,2i-1),\wedge^2,k})=((2n)^{m-2k+1},n^{4k-2}).
\end{equation}
We conclude from \eqref{5.2} that $(2n,1^{2mn})\circ ((2n)^{m-2k+1},n^{4k-2})$ supports $\mathcal{E}^{\Sp^{(2)}_{4ni}}_{\Delta(\tau,2i),\wedge^2,k}$. We used the notion of composition of unipotent classes, defined in the end of Sec. 1 in \cite{GRS03}. Note that the Fourier coefficient used to define $\mathcal{D}^\phi_{\psi,n(2i-1)}$ in \eqref{5.1'} corresponds to the partition $(2n,1^{2n(2i-1)})$ of $4ni$.
By \cite{GRS03}, Lemma 6,
$((2n)^{m-2k+2},n^{4k-2})$ supports $\mathcal{E}^{\Sp^{(2)}_{4ni}}_{\Delta(\tau,2i),\wedge^2,k}$. This proves that, for all $1\leq k\leq i$,
\begin{equation}\label{5.3}
\mathcal{O}(\mathcal{E}^{\Sp^{(2)}_{4ni}}_{\Delta(\tau,2i),\wedge^2,k})=((2n)^{2i-2k+1},n^{4k-2}).
\end{equation}
Similarly, from \eqref{5.1''} and \eqref{5.3}, it follows that $(2n,1^{4ni})\circ ((2n)^{2i-2k+1},n^{4k-2})$ supports $\mathcal{E}^{\Sp_{2n(2i+1)}}_{\Delta(\tau,2i+1),\wedge^2,k}$. By \cite{GRS03}, Lemma 6,
$((2n)^{2i-2k+2},n^{4k-2})$ supports $\mathcal{E}^{\Sp_{2n(2i+1)}}_{\Delta(\tau,2i+1),\wedge^2,k}$, and hence, for all $1\leq k\leq i=\frac{m+1}{2}$,
\begin{equation}\label{5.5'}
\mathcal{O}(\mathcal{E}^{\Sp_{2n(2i+1)}}_{\Delta(\tau,2i+1),\wedge^2,k})=((2n)^{(2i+1)-2k+1},n^{4k-2}).
\end{equation}
For $k=i+1=\frac{m+3}{2}$, we know from Theorem \ref{thm 4.1} that
$$
\mathcal{O}(\mathcal{E}^{\Sp_{2n(m+2)}}_{\Delta(\tau,m+2),\wedge^2,\frac{m+3}{2}})=(n^{2(m+2)}).
$$
Thus, \eqref{5.5'} is valid for $k=\frac{m+3}{2}$, as well. This proves parts I.2, I.3 of the theorem.

We now turn to the poles of $E^{\Sp^{(2)}_{2n(2i+1)}}(f_{\Delta(\tau,2i+1)\gamma_\psi,s})$.
By Theorem \ref{thm 7.1}, for all $h\in \Sp_{4ni}(\BA)$,	
\begin{equation}\label{5.2'}
\mathcal{D}^\phi_{\psi,2ni}(E^{\Sp^{(2)}_{2n(2i+1)}}(f_{\Delta(\tau,2i+1)\gamma_\psi,s}))(h)=E^{\Sp_{4ni}}(\Lambda(f_{\Delta(\tau,2i+1)\gamma_\psi,s},\phi))(h).
\end{equation}
Now, we argue exactly as we did right after \eqref{5.1'}, and get from \eqref{5.2'} that\\
 $E^{\Sp^{(2)}_{2n(2i+1)}}(f_{\Delta(\tau,2i+1)\gamma_\psi,s})$ has a pole at $s=j$, for all $1\leq j\leq i$ (as the section varies).
Next, by Theorem \ref{thm 7.1}, for all $h\in \Sp^{(2)}_{2n(2i+1)}(\BA)$,
\begin{equation}\label{5.2''}
\mathcal{D}^\phi_{\psi,2n(i+1)}(E^{\Sp_{4n(i+1)}}(f_{\Delta(\tau,2i+2),s}))(h)=E^{\Sp^{(2)}_{2n(2i+1)}}(\Lambda(f_{\Delta(\tau,2i+2),s},\phi))(h).
\end{equation}
Again, we argue exactly as we did right after \eqref{5.1''}, and get from \eqref{5.2''} that each pole $s=j$, $1\leq j\leq i$ of $E^{\Sp^{(2)}_{2n(2i+1)}}(f_{\Delta(\tau,2i+1)\gamma_\psi,s})$ is simple. Finally, let $s_0$, with $Re(s_0)\geq 0$ be a pole of $E^{\Sp^{(2)}_{2n(2i+1)}}(f'_{\Delta(\tau,2i+1)\gamma_\psi,s})$. Then we conclude from \eqref{5.2''} that it is a pole of $E^{\Sp_{4n(i+1)}}(f_{\Delta(\tau,2i+2),s})$ as the section varies, and hence there is a $1\leq j\leq i+1$, such that $s_0=j$. As before, we must have $1\leq j\leq i$, since $s=i+1$ is not a pole of $E^{\Sp^{(2)}_{2n(2i+1)}}(f_{\Delta(\tau,2i+1)\gamma_\psi,s})$. This can be seen by the proof of Prop. \ref{prop (0)}.

We now prove Parts I.1, I.4. The proof is similar to the last case.
We prove Part I.1 by induction on $m$ (even) and Part I.3 follows from the proof. When $m=2$, Part I.1 is a special case of Theorem \ref{thm 4.1}. Assume by induction that, for $m=2i$ even, and $1\leq k\leq \frac{m}{2}=i$,
\begin{equation}\label{5.6.1}
\mathcal{O}(\mathcal{E}^{\Sp_{2mn}}_{\Delta(\tau,m),\wedge^2,k})=((2n)^{m-2k},n^{4k}).
\end{equation}
From \eqref{5.6.1}, we get that $(2n,1^{2mn})\circ ((2n)^{m-2k},n^{4k})$ supports $\mathcal{E}^{\Sp^{(2)}_{2n(2i+1)}}_{\Delta(\tau,2i+1),\wedge^2,k}$.
By \cite{GRS03}, Lemma 6,
$((2n)^{m-2k+1},n^{4k})$ supports $\mathcal{E}^{\Sp^{(2)}_{2n(2i+1)}}_{\Delta(\tau,2i+1),\wedge^2,k}$. This proves that, for all $1\leq k\leq i$,
\begin{equation}\label{5.8}
\mathcal{E}^{\Sp^{(2)}_{2n(2i+1)}}_{\Delta(\tau,2i+1),\wedge^2,k}=((2n)^{2i+1-2k},n^{4k}).
\end{equation}
Similarly, from \eqref{5.2''} and \eqref{5.8}, it follows that $(2n,1^{2n(2i+1)})\circ ((2n)^{2i+1-2k},n^{4k})$ supports $\mathcal{E}^{\Sp_{4n(i+1)}}_{\Delta(\tau,2i+2),\wedge^2,k}$. By \cite{GRS03}, Lemma 6,
$((2n)^{2(i+1)-2k},n^{4k})$ supports\\ $\mathcal{E}^{\Sp_{4n(i+1)}}_{\Delta(\tau,2i+2),\wedge^2,k}$, and hence, for all $1\leq k\leq i=\frac{m+2}{2}$,
\begin{equation}\label{5.4'}
\mathcal{E}^{\Sp_{4n(i+1)}}_{\Delta(\tau,2i+2),\wedge^2,k}=((2n)^{2(i+1)-2k},n^{4k}).
\end{equation}
For $k=i+1=\frac{m+2}{2}$, we know from Theorem \ref{thm 4.1} that
$$
\mathcal{O}(\mathcal{E}^{\Sp_{2n(m+2)}}_{\Delta(\tau,m+2),\wedge^2,\frac{m+2}{2}})=(n^{2(m+2)}).
$$
Thus, \eqref{5.4'} is valid for $k=\frac{m+2}{2}$, as well. This proves parts I.1, I.4 of the theorem.\\
\vspace{0.1cm}

II. $L(\tau,\vee^2,s)$ has a pole at $s=1$, and $\omega_\tau=1$.\\

Assume that $L(\tau,\vee^2,s)$ has a pole at $s=1$, and $\omega_\tau=1$. The statement of the theorem in detail is the following:

The set of points \eqref{1.4}, $e^{(2)}_{k,m}(\vee^2)$,
is the set of poles of $E(f_{\Delta(\tau,m)\gamma_\psi,s})$, as the section varies, in $Re(s)\geq 0$, and they are all simple. (We know from \cite{JLZ13} that the set of points \eqref{1.3}, $e_{k,m}(\vee^2)$, is the set of poles of  $E(f_{\Delta(\tau,m),s})$, as the section varies, in $R(s)\geq 0$, and they are all simple.)
We have

\begin{enumerate}
	
\item For $m$ even and $1\leq k\leq \frac{m}{2}$ ($e^{\Sp_{2mn}}_{k,m}(\vee^2)=k-\frac{1}{2}$)
$$
\mathcal{O}(\mathcal{E}^{\Sp_{2mn}}_{\Delta(\tau,m),\vee^2,k})=
\begin{cases}	
((2n)^{m-2k+1},n^{4k-2}),\ \ \ \ \ \ \ \ \ \ \ \ \ n\ even\\
((2n)^{m-2k+1},(n+1)^{2k-1},(n-1)^{2k-1}),\ \ n\ odd
\end{cases}
$$

\item 	For $m$ odd and $1\leq k\leq \frac{m-1}{2}$ ($e^{\Sp_{2mn}}_{k,m}(\vee^2)=k$)
$$
\mathcal{O}(\mathcal{E}^{\Sp_{2mn}}_{\Delta(\tau,m),\vee^2,k})=
\begin{cases}
((2n)^{m-2k},n^{4k}),\ \ \ \ \ \ \ \ \ \ \ \ \ \ n\ even\\
((2n)^{m-2k},(n+1)^{2k},(n-1)^{2k}),\ \ n\ odd
\end{cases}
$$

\item 	For $m$ even and $1\leq k\leq \frac{m}{2}$ ($e^{\Sp^{(2)}_{2mn}}_{k,m}(\vee^2)=k$)
$$
\mathcal{O}(\mathcal{E}^{\Sp^{(2)}_{2mn}}_{\Delta(\tau,m)\gamma_\psi,\vee^2,k})=
\begin{cases}
((2n)^{m-2k},n^{4k}),\ \ \ \ \ \ \ \ \ \ \ \ \ \ n\ even\\
((2n)^{m-2k},(n+1)^{2k},(n-1)^{2k}),\ \ n\ odd
\end{cases}
$$

\item For $m$ odd and $1\leq k\leq \frac{m+1}{2}$ ($e^{\Sp^{(2)}_{2mn}}_{k,m}(\vee^2)=k-\frac{1}{2}$)
$$
\mathcal{O}(\mathcal{E}^{\Sp^{(2)}_{2mn}}_{\Delta(\tau,m)\gamma_\psi,\vee^2,k})=
\begin{cases}
((2n)^{m-2k+1},n^{4k-2}),\ \ \ \ \ \ \ \ \ \ \ \ \ \ n\ even\\
((2n)^{m-2k+1},(n+1)^{2k-1},(n-1)^{2k-1}),\ \ n\ odd
\end{cases}
$$	
	
\end{enumerate}

We first prove the statement about the poles of $E^{\Sp^{(2)}_{2n(2i-1)}}(f_{\Delta(\tau,2i-1)\gamma_\psi,s})$. The proof is similar to the one in the previous case. By Theorem \ref{thm 7.1}, for all $h\in \Sp_{4n(i-1)}(\BA)$,	
\begin{equation}\label{5.3''}
\mathcal{D}^\phi_{\psi,2n(i-1)}(E^{\Sp^{(2)}_{2n(2i-1)}}(f_{\Delta(\tau,2i-1)\gamma_\psi,s}))(h)=E^{\Sp_{4n(i-1)}}(\Lambda(f_{\Delta(\tau,2i-1)\gamma_\psi,s},\phi))(h).
\end{equation}
Now, we argue exactly as we did right after \eqref{5.1'}, \eqref{5.2'}, and get from \eqref{5.3''} that $E^{\Sp^{(2)}_{2n(2i-1)}}(f_{\Delta(\tau,2i-1)\gamma_\psi,s})$ has a pole at $s=\frac{2j-1}{2}$, for all $1\leq j\leq i-1$. In case $i=1$, the l.h.s. of \eqref{5.3''} is the Whittaker coefficient of $E^{\Sp^{(2)}_{2n}}(f_{\tau\gamma_\psi,s})$, which is holomorphic. By Theorem \ref{thm 7.1}, for all $h\in \Sp^{(2)}_{2n(2i-1)}(\BA)$,	
\begin{equation}\label{5.3'}
\mathcal{D}^\phi_{\psi,n(2i-1)}(E^{\Sp_{4ni}}(f_{\Delta(\tau,2i),s}))(h)=E^{\Sp^{(2)}_{2n(2i-1)}}(\Lambda(f_{\Delta(\tau,2i),s},\phi))(h).
\end{equation}
As before, since the points $s=\frac{2j-1}{2}$, for all $1\leq j\leq i$ are simple poles of $E^{\Sp_{4ni}}(f_{\Delta(\tau,2i)),s})$, we conclude that each point $s=\frac{2j-1}{2}$, $1\leq j\leq i-1$, is a simple pole of $E^{\Sp^{(2)}_{2n(2i-1)}}(f_{\Delta(\tau,2i-1)\gamma_\psi,s})$, as the section varies. Note that by Prop. \ref{prop (0)}, $s=\frac{2i-1}{2}$ is also a simple pole of the last series. Let $s_0$, with $Re(s_0)\geq 0$ be a pole of $E^{\Sp^{(2)}_{2n(2i-1)}}(f'_{\Delta(\tau,2i-1)\gamma_\psi,s})$. From \eqref{5.3'}, we get that $s_0=\frac{2j-1}{2}$, for some $1\leq j\leq i$. Note that in the case $i=1$, we get that $E^{\Sp^{(2)}_{2n}}(f_{\tau\gamma_\psi,s})$ has a single pole in $Re(s)\geq 0$, which is $s=\frac{1}{2}$.

We will now prove Part II.4 of the theorem by induction on $m$. The proof then implies part II.1. The proof follows the same lines of the last proof. We note again that the indicated partition in Parts II.1, II.4 bounds $\mathcal{O}(\mathcal{E}^{\Sp^{(\epsilon)}_{2mn}}_{\Delta(\tau,m)\gamma^{(\epsilon)}_\psi,\vee^2,k})$.
This is proved in \cite{GS20}, Prop. 3.2, when $m$ is even. The proof when $m$ is odd is entirely the same. Thus, we need to show that the residual Eisenstein series in these two cases admit nontrivial Fourier coefficients corresponding to the indicated top partition.

We start with $E^{\Sp^{(2)}_{2n}}(f_{\tau,s})$. By Prop. \ref{prop (0)}, it has a simple pole at $s=\frac{1}{2}$, and by Theorem \ref{thm 4.2}(1),
$$
 \mathcal{O}(\mathcal{E}^{\Sp^{(2)}_{2n}}_{\tau\gamma_\psi,\vee^2,1})=
 \begin{cases}
(n^2),\ \ \ \ \ \ \ \ n\ even\\
(n+1,n-1),\ \ n\ odd
\end{cases}
$$
This proves Part II.4 of the theorem for $m=1$. Assume by induction that, for $m=2i-1$ odd, and $1\leq k\leq \frac{m+1}{2}=i$,
$$
 \mathcal{O}(\mathcal{E}^{\Sp^{(2)}_{2n(2i-1)}}_{\Delta(\tau,2i-1)\gamma_\psi,\vee^2,k})=
 \begin{cases}
 ((2n)^{m-2k+1},n^{4k-2}),\ \ \ \ \ \ \ \ \ \ \ \ n\ even\\
((2n)^{m-2k+1},(n+1)^{2k-1},(n-1)^{2k-1}),\ \ n\ odd
\end{cases}
$$
Let us write the relation \eqref{5.3''} with $i+1$ instead of $i$, that is for all $h\in \Sp_{4ni}(\BA)$,	
\begin{equation}\label{5.4''}
\mathcal{D}^\phi_{\psi,2ni}(E^{\Sp^{(2)}_{2n(2i+1)}}(f_{\Delta(\tau,2i+1)\gamma_\psi,s}))(h)=E^{\Sp_{4ni}}(\Lambda(f_{\Delta(\tau,2i+1)\gamma_\psi,s},\phi))(h).
\end{equation}

We conclude from \eqref{5.3'} and the induction assumption, always using Cor. \ref{cor (1')} and Prop. \ref{prop (2)}, that $(2n,1^{2mn})\circ ((2n)^{m-2k+1},n^{4k-2})$, when $n$ is even, and $(2n,1^{2mn})\circ ((2n)^{m-2k+1},(n+1)^{2k-1}, (n-1)^{2k-1})$, when $n$ is odd, support $\mathcal{E}^{\Sp_{4ni}}_{\Delta(\tau,2i),\vee^2,k}$.  By \cite{GRS03}, Lemma 6, $((2n)^{2i-2k+1},n^{4k-2})$, when $n$ is even, and $ ((2n)^{2i-2k+1},(n+1)^{2k-1}, (n-1)^{2k-1})$, when $n$ is odd, support $\mathcal{E}^{\Sp_{4ni}}_{\Delta(\tau,2i),\vee^2,k}$.
Hence, for all $1\leq k\leq \frac{2i}{2}=i$,
$$
\mathcal{O}(\mathcal{E}^{\Sp_{4ni}}_{\Delta(\tau,2i),\vee^2,k})=
\begin{cases}
((2n)^{2i-2k+1},n^{4k-2}),\ \ \ \ \ \ \ \ \ \ \ \ n\ even\\
((2n)^{2i-2k+1},(n+1)^{2k-1},(n-1)^{2k-1}),\ \ n\ odd
\end{cases}
$$
Now we conclude from the last equality and \eqref{5.4''} (using Cor. \ref{cor (1')} and Prop. \ref{prop (2)}) that $(2n,1^{4ni})\circ ((2n)^{2i-2k+1},n^{4k-2})$, when $n$ is even, and $(2n,1^{4ni})\circ ((2n)^{2i-2k+1},(n+1)^{2k-1},(n-1)^{2k-1})$, when $n$ is odd, support $\mathcal{E}^{\Sp^{(2)}_{2n(2i+1)}}_{\Delta(\tau,2i+1)\gamma_\psi,\vee^2,k}$. By \cite{GRS03}, Lemma 6,  $((2n)^{2i-2k+2},n^{4k-2})$, when $n$ is even, and $ ((2n)^{2i-2k+2},(n+1)^{2k-1},(n-1)^{2k-1})$, when $n$ is odd, support $\mathcal{E}^{\Sp^{(2)}_{2n(2i+1)}}_{\Delta(\tau,2i+1)\gamma_\psi,\vee^2,k}$. Hence, for all $1\leq k<\frac{(2i+1)+1}{2}=i+1$,
$$
\mathcal{O}(\mathcal{E}^{\Sp^{(2)}_{2n(2i+1)}}_{\Delta(\tau,2i+1)\gamma_\psi,\vee^2,k})=
\begin{cases}
((2n)^{(2i+1)-2k+1},n^{4k-2}),\ \ \ \ \ \ \ \ \ \ \ \ n\ even\\
((2n)^{(2i+1)-2k+1},(n+1)^{2k-1},(n-1)^{2k-1}),\ \ n\ odd
\end{cases}
$$
The last equality is true for $k=\frac{(2i+1)+1}{2}=i+1$, as well, by Theorem \ref{thm 4.2}(1). This proves parts II.4, II.1 of the theorem.

The proof of the determination of poles of $E^{\Sp^{(2)}_{4ni}}(f_{\Delta(\tau,2i)\gamma_\psi,s})$ and the proof of Parts II.3, II.2 is entirely similar. One proves Part II.3 by induction on $m$ (even) and part II.2 follows from the proof. Note that part II.3 for $m=2$ is a special case of Theorem \ref{thm 4.2}. We omit the details.\\
\vspace{0.1cm}

III. $L(\tau,\vee^2,s)$ has a pole at $s=1$, and $\omega_\tau\neq 1$.\\

Assume that $L(\tau,\vee^2,s)$ has a pole at $s=1$, and $\omega_\tau\neq 1$.
Note that the proof in the last part that each point $e^{(2)}_{k,m}(\vee^2)$
is a simple pole of $E(f_{\Delta(\tau,m)\gamma_\psi,s})$, as the section varies did not depend on $\omega_\tau$ being trivial or not.
The statement of the theorem in detail is the following:

\begin{enumerate}
	
	\item For $m$ even and $1\leq k\leq \frac{m}{2}$ ($e^{\Sp_{2mn}}_{k,m}(\vee^2)=k-\frac{1}{2}$)
	$$
	\mathcal{O}(\mathcal{E}^{\Sp_{2mn}}_{\Delta(\tau,m),\vee^2,k})=
	\begin{cases}	
	((2n)^{m-2k+1},(n+2)^{2k-1}, (n-2)^{2k-1}),\ \ n\ even\\
	((2n)^{m-2k+1},(n+1)^{2k-1},(n-1)^{2k-1}),\ \ n\ odd
	\end{cases}
	$$

	\item 	For $m$ odd and $1\leq k\leq \frac{m-1}{2}$ ($e^{\Sp_{2mn}}_{k,m}(\vee^2)=k$)
	$$
	\mathcal{O}(\mathcal{E}^{\Sp_{2mn}}_{\Delta(\tau,m),\vee^2,k})=
	\begin{cases}
	((2n)^{m-2k},(n+2)^{2k}, (n-2)^{2k}),\ \ n\ even\\
	((2n)^{m-2k},(n+1)^{2k},(n-1)^{2k}),\ \ n\ odd
	\end{cases}
	$$

	\item 	For $m$ even and $1\leq k\leq \frac{m}{2}$ ($e^{\Sp^{(2)}_{2mn}}_{k,m}(\vee^2)=k$)
	$$
	\mathcal{O}(\mathcal{E}^{\Sp^{(2)}_{2mn}}_{\Delta(\tau,m)\gamma_\psi,\vee^2,k})=
	\begin{cases}
((2n)^{m-2k},(n+2)^{2k}, (n-2)^{2k}),\ \ n\ even\\
	((2n)^{m-2k},(n+1)^{2k},(n-1)^{2k}),\ \ n\ odd
	\end{cases}
	$$

	\item For $m$ odd and $1\leq k\leq \frac{m+1}{2}$ ($e^{\Sp^{(2)}_{2mn}}_{k,m}(\vee^2)=k-\frac{1}{2}$)
	$$
	\mathcal{O}(\mathcal{E}^{\Sp^{(2)}_{2mn}}_{\Delta(\tau,m)\gamma_\psi,\vee^2,k})=
	\begin{cases}
	((2n)^{m-2k+1},(n+2)^{2k-1}, (n-2)^{2k-1}),\ \ n\ even\\
	((2n)^{m-2k+1},(n+1)^{2k-1},(n-1)^{2k-1}),\ \ n\ odd
	\end{cases}
	$$	
	
\end{enumerate}

When $n$ is odd, this part is the same as the last part. Thus, assume that $n=2n'$ is even.
In each case, the indicated partition bounds $\mathcal{O}(\mathcal{E}^{\Sp^{(\epsilon)}_{2mn}}_{\Delta(\tau,m)\gamma^{(\epsilon)}_\psi,\vee^2,k})$. The proof is similar to that of Prop. 3.2 in \cite{GS20}. Let us sketch it for $\Sp^{(2)}_{2mn}$. The case of the residue at $s=\frac{m}{2}$ was proved in the beginning of the proof of Theorem \ref{thm 4.2}(2). The proof for $1\leq k<[\frac{m+1}{2}]$ is similar. Let $v$ be a finite place, where $\tau_v$ is unramified, and its central character is the unique, unramified, nontrivial character $\lambda_v$. Write $\tau_v$ as in \eqref{4.19.1}.
Then the unramified constituent of the factor at $v$ of 	
$\mathcal{E}^{\Sp^{(2)}_{2nm}}_{\Delta(\tau,m)\gamma_\psi,\vee^2,k}$  is the unramified constituent of the following parabolic inductions, according to whether $m$ is even, or odd. When $m$ is even, this is the parabolic induction from the following character of $Q^{(2)}_{(m+2k)^{n'-1},(m-2k)^{n'-1},m^2}(F_v)$,
\begin{equation}\label{5.6}
[\otimes_{i=1}^{n'-1}(\chi_i\circ det_{\GL_{m+2k}})\otimes (\chi_i\circ det_{\GL_{m-2k}}) \otimes (\lambda_v\circ det_{\GL_m})|det_{\GL_m}|^k\otimes |det_{\GL_m}|^k]\gamma_{\psi_v};
\end{equation}
When $m$ is odd, this is the parabolic induction from the following character of $Q^{(2)}_{(m+2k-1)^{n'-1},(m-2k+1)^{n'-1},m^2}(F_v)$,
\begin{multline}\label{5.7}
[\otimes_{i=1}^{n'-1}(\chi_i\circ det_{\GL_{m+2k-1}})\otimes (\chi_i\circ det_{\GL_{m-2k+1}}) \otimes (\lambda_v\circ det_{\GL_m})|det_{\GL_m}|^{k-\frac{1}{2}}\otimes\\ \otimes|det_{\GL_m}|^{k-\frac{1}{2}}]\gamma_{\psi_v};
\end{multline}

As in the proof of Prop. 3.2 in \cite{GS20}, we conclude, from \eqref{5.6}, \eqref{5.7}, that all symplectic partitions of $2mn$ corresponding to nontrivial Fourier coefficients on  $\mathcal{E}^{\Sp^{(2)}_{2nm}}_{\Delta(\tau,m)\gamma_\psi,\vee^2,k}$ are bounded by the induced nilpotent orbit corresponding to \eqref{5.6}, \eqref{5.7}, and this  corresponds to the partition $((2n)^{m-2k}, (n+2)^{2k}, (n-2)^{2k}))$, when $m$ is even, and $((2n)^{m-2k+1}, (n+2)^{2k-1}, (n-2)^{2k-1}))$, when $m$ is odd. See \cite{CM93}, Chapter 7. It remains to show that the residual Eisenstein series admits in each case a nontrivial Fourier coefficient corresponding to the indicated top partition.

We now prove Part III.4 by induction on $m$ (even) and Part III.1 follows from the proof. The proof is the same as in the last part, except that we need to replace at each place in the proof (when $n$ is even) the partition $((2n)^{m-2k+1},n^{4k-2})$ by the partition $((2n)^{m-2k+1},(n+2)^{2k-1}, (n-2)^{2k-1})$. We start with $E^{\Sp^{(2)}_{2n}}(f_{\tau,s})$. It has a pole at $s=\frac{1}{2}$, and by Theorem \ref{thm 4.2}(2),
$$
\mathcal{O}(\mathcal{E}^{\Sp^{(2)}_{2n}}_{\tau\gamma_\psi,\vee^2,1})=(n+2,n-2).
$$
This proves Part III.4 of the theorem for $m=1$. Assume by induction that, for $m=2i-1$ odd, and $1\leq k\leq \frac{m+1}{2}=i$ (recall that $n$ is even),
$$
\mathcal{O}(\mathcal{E}^{\Sp^{(2)}_{2n(2i-1)}}_{\Delta(\tau,2i-1)\gamma_\psi,\vee^2,k})=((2n)^{m-2k+1},(n+2)^{2k-1}, (n-2)^{2k-1}).
$$
As in the last part, we conclude from \eqref{5.3'} and the induction assumption, that $(2n,1^{2mn})\circ ((2n)^{m-2k+1},(n+2)^{2k-1}, (n-2)^{2k-1})$ supports $\mathcal{E}^{\Sp_{4ni}}_{\Delta(\tau,2i),\vee^2,k}$.  By \cite{GRS03}, Lemma 6, $((2n)^{2i-2k+1},(n+2)^{2k-1}, (n-2)^{2k-1})$ supports $\mathcal{E}^{\Sp_{4ni}}_{\Delta(\tau,2i),\vee^2,k}$.
Hence, for all $1\leq k\leq \frac{2i}{2}=i$,
$$
\mathcal{O}(\mathcal{E}^{\Sp_{4ni}}_{\Delta(\tau,2i),\vee^2,k})=((2n)^{2i-2k+1},(n+2)^{2k-1}, (n-2)^{2k-1}).
$$
By Theorem \ref{thm 7.1}, we have the relation \eqref{5.3''}. Again, we conclude from the last equality and \eqref{5.3''} that $(2n,1^{4ni})\circ ((2n)^{2i-2k+1},(n+2)^{2k-1}, (n-2)^{2k-1})$, supports $\mathcal{E}^{\Sp^{(2)}_{2n(2i+1)}}_{\Delta(\tau,2i+1)\gamma_\psi,\vee^2,k}$. By \cite{GRS03}, Lemma 6,  $((2n)^{2i-2k+2},(n+2)^{2k-1}, (n-2)^{2k-1})$ supports $\mathcal{E}^{\Sp^{(2)}_{2n(2i+1)}}_{\Delta(\tau,2i+1)\gamma_\psi,\vee^2,k}$. Hence, for all $1\leq k<\frac{(2i+1)+1}{2}=i+1$,
$$
\mathcal{O}(\mathcal{E}^{\Sp^{(2)}_{2n(2i+1)}}_{\Delta(\tau,2i+1)\gamma_\psi,\vee^2,k})=
((2n)^{(2i+1)-2k+1},(n+2)^{2k-1}, (n-2)^{2k-1}).
$$
The last equality is true for $k=\frac{(2i+1)+1}{2}=i+1$, as well, by Theorem \ref{thm 4.2}(2). This proves Parts III.4, III.1 of the theorem. We leave the proof of Parts III.3, III.2 to the reader.

\end{document}